\documentclass[11pt]{amsart}
\usepackage{amsmath,amssymb,wasysym}
\newtheorem{theorem}{Theorem}[section]
\newtheorem{proposition}[theorem]{Proposition}
\newtheorem{corollary}[theorem]{Corollary}
\newtheorem{lemma}[theorem]{Lemma}
\theoremstyle{definition}
\newtheorem{definition}[theorem]{Definition}
\newtheorem{remark}[theorem]{Remark}
\newtheorem{assumption}[theorem]{Assumption}

\begin{document}

\title[Gromov's rigidity theorem]{On Gromov's rigidity theorem for polytopes with acute angles}
\author{Simon Brendle and Yipeng Wang}
\address{Columbia University \\ 2990 Broadway \\ New York NY 10027 \\ USA}
\address{Columbia University \\ 2990 Broadway \\ New York NY 10027 \\ USA}
\begin{abstract}
In his ``Four Lectures", Gromov conjectured a scalar curvature extremality property of convex polytopes. Moreover, Gromov outlined a proof of the conjecture in the special case when the dihedral angles are acute. Gromov's argument relies on Dirac operator techniques together with a smoothing construction. In this paper, we give the details of such a smoothing construction, thereby providing a detailed proof of Gromov's theorem. 
\end{abstract}
\thanks{The author was supported by the National Science Foundation under grants DMS-2103573 and DMS-2403981 and by the Simons Foundation. He acknowledges the hospitality of T\"ubingen University, where part of this work was carried out.}
\maketitle

\section{Introduction}

In this paper, we prove a scalar curvature rigidity property for polytopes with acute angles. Suppose that $n \geq 3$ is an integer, and $\Omega$ is a compact, convex polytope in $\mathbb{R}^n$ with non-empty interior. We write $\Omega = \bigcap_{m=0}^q \{u_m \leq 0\}$ where $u_0,\hdots,u_q$ are non-constant linear functions defined on $\mathbb{R}^n$. After eliminating redundant inequalities, we may assume that the following condition is satisfied.

\begin{assumption}
\label{no.redundant.inequalities}
For each $0 \leq k \leq q$, the set 
\[\{u_k > 0\} \cap \bigcap_{m \in \{0,1,\hdots,q\} \setminus \{k\}} \{u_m \leq 0\}\] 
is non-empty. 
\end{assumption}

Without loss of generality, we assume that the following condition is satisfied.

\begin{assumption}
\label{gradient.of.u_k}
For each $0 \leq k \leq q$, the gradient of the function $u_k$ with respect to the Euclidean metric is given by a unit vector $N_k \in S^{n-1}$.
\end{assumption}

\begin{theorem}[cf. M.~Gromov \cite{Gromov2}, Section 3.18]
\label{main.thm}
Suppose that $n \geq 3$ is an integer, and $\Omega = \bigcap_{m=0}^q \{u_m \leq 0\}$ is a compact, convex polytope in $\mathbb{R}^n$ with non-empty interior. Suppose that Assumptions \ref{no.redundant.inequalities} and \ref{gradient.of.u_k} are satisfied. Let $g$ be a Riemannian metric on $\mathbb{R}^n$. For each $k \in \{0,1,\hdots,q\}$, we denote by $\nu_k = \frac{\nabla u_k}{|\nabla u_k|}$ the unit normal vector to the level sets of $u_k$ with respect to the metric $g$. We assume that the following conditions are satisfied:
\begin{itemize}
\item The scalar curvature of $g$ is nonnegative at each point in $\Omega$.
\item For each $k \in \{0,1,\hdots,q\}$, the mean curvature of the hypersurface $\{u_k=0\}$ with respect to $g$ is nonnegative at each point in $\Omega \cap \{u_k=0\}$. 
\item If $j$ and $k$ are integers with $0 \leq j < k \leq q$ and $x$ is a point in $\Omega$ with $u_j(x) = u_k(x) = 0$, then $\langle \nu_j,\nu_k \rangle \leq \langle N_j,N_k \rangle \leq 0$ at the point $x$. 
\end{itemize}
Then:
\begin{itemize}
\item The Riemann curvature tensor of $g$ vanishes at each point in $\Omega$.
\item For each $k \in \{0,1,\hdots,q\}$, the second fundamental form of the hypersurface $\{u_k=0\}$ with respect to $g$ vanishes at each point in $\Omega \cap \{u_k=0\}$. 
\item If $j$ and $k$ are integers with $0 \leq j < k \leq q$ and $x$ is a point in $\Omega$ with $u_j(x) = u_k(x) = 0$, then $\langle \nu_j,\nu_k \rangle = \langle N_j,N_k \rangle$ at the point $x$. 
\end{itemize}
\end{theorem}

Theorem \ref{main.thm} is a special case of Gromov's Dihedral Rigidity Conjecture. Theorem \ref{main.thm} was stated by Gromov in his ``Four Lectures'' (see \cite{Gromov2}, Section 3.18). Moreover, Gromov sketched a proof based on Dirac operator techniques and a smoothing construction. The main purpose of this paper is to give the details of the smoothing construction. 

In Section \ref{smoothing.of.polytope}, we explain how a given polytope $\Omega$ may be approximated by a convex domain $\hat{\Omega} \subset \Omega$ with smooth boundary $\partial \hat{\Omega} = \hat{\Sigma}$. 

In Section \ref{construction.hat.N}, we construct a smooth map $\hat{N}: \hat{\Sigma} \to S^{n-1}$ which is homotopic to the Gauss map of $\hat{\Sigma}$. 

In Section \ref{controlling.angles}, we explain how to control certain angles under the smoothing procedure (see Lemma \ref{preservation.of.angle.inequality.under.smoothing} and Proposition \ref{approximate.angle.inequality.3} below). Here, we use in an essential way the assumption that $\Omega$ has acute angles.

In Section \ref{controlling.mean.curvature}, we prove pointwise estimates for $\max \{\|d\hat{N}\|_{\text{\rm tr}} - H,0\}$ in various subsets of $\hat{\Sigma}$, where $H$ denotes the mean curvature of $\hat{\Sigma}$ with respect to the metric $g$. From this, we deduce that a suitable Morrey norm of $\max \{\|d\hat{N}\|_{\text{\rm tr}} - H,0\}$ can be made arbitrarily small (see Corollary \ref{Morrey.norm} below). 

In Section \ref{proof.of.main.thm}, we complete the proof of Theorem \ref{main.thm}. The proof follows the arguments in \cite{Brendle}. It uses the bound for the Morrey norm of $\max \{\|d\hat{N}\|_{\text{\rm tr}} - H,0\}$ established in Section \ref{controlling.mean.curvature}, together with a deep estimate due to Fefferman and Phong \cite{Fefferman-Phong}.

We refer to \cite{Brendle}, \cite{Gromov1}, \cite{Gromov2}, \cite{Gromov3}, \cite{Li1}, \cite{Li2} for related work on Gromov's Dihedral Rigidity Conjecture. We note that Wang, Xie, and Yu \cite{Wang-Xie-Yu} have proposed a different approach to this problem which is based on the study of Dirac operators on manifolds with corners; their work is currently in the process of verification.

\section{Smoothing convex polytopes with acute angles}

\label{smoothing.of.polytope}

Throughout this section, we assume that $n \geq 3$ is an integer, and $\Omega$ is a compact, convex polytope in $\mathbb{R}^n$ with non-empty interior. As above, we write $\Omega = \bigcap_{m=0}^q \{u_m \leq 0\}$ where $u_0,\hdots,u_q$ are non-constant linear functions defined on $\mathbb{R}^n$. Throughout this section, we assume that Assumption \ref{no.redundant.inequalities} and \ref{gradient.of.u_k} are satisfied. Moreover, we assume that the following assumption is satisfied.

\begin{assumption}
\label{angles.bounded.by.pi/2}
Let $0 \leq j < k \leq q$. If there exists a point $x \in \Omega$ satisfying $u_j(x)=u_k(x)=0$, then $\langle N_j,N_k \rangle \leq 0$. 
\end{assumption}

\begin{lemma}
\label{approximate.angle.inequality.1}
Given $\varepsilon \in (0,1)$, we can find a small positive real number $\delta$ (depending on $\varepsilon$) with the following property. Suppose that $0 \leq j < k \leq q$. If there exists a point $x \in \Omega$ with $-\delta \leq u_j(x) \leq 0$ and $-\delta \leq u_k(x) \leq 0$, then $\langle N_j,N_k \rangle \leq \varepsilon$.
\end{lemma}

\textbf{Proof.} 
Suppose that the assertion is false. We consider a sequence of counterexamples and pass to the limit. In the limit, we obtain a pair of integers $0 \leq j < k \leq q$ and a point $x \in \Omega$ such that $u_j(x)=u_k(x)=0$ and $\langle N_j,N_k \rangle \geq \varepsilon$. This contradicts Assumption \ref{angles.bounded.by.pi/2}. This completes the proof of Lemma \ref{approximate.angle.inequality.1}. \\

\begin{lemma}
\label{Lambda}
We can find a large constant $\Lambda>1$ with the following property. Suppose that $x$ is a point in $\Omega$ and $a_0,\hdots,a_q$ are nonnegative real numbers such that $a_i = 0$ for all $0 \leq i \leq q$ satisfying $u_i(x) \leq -2\Lambda^{-1}$. Then $\sum_{i=0}^q a_i \leq \Lambda \, \big | \sum_{i=0}^q a_i \, N_i \big |$. 
\end{lemma} 

\textbf{Proof.} 
Suppose that the assertion is false. We consider a sequence of counterexamples and pass to the limit. In the limit, we obtain a point $x \in \Omega$ and a collection of nonnegative real numbers $a_0,\hdots,a_q$ with the following properties: 
\begin{itemize} 
\item $\sum_{i=0}^q a_i = 1$.
\item $\sum_{i=0}^q a_i \, N_i = 0$. 
\item $a_i = 0$ for all $0 \leq i \leq q$ satisfying $u_i(x) < 0$. 
\end{itemize}
For abbreviation, let $I := \{0,1,\hdots,q\}$ and $I_0 := \{i \in I: u_i(x)=0\}$. Since $\Omega$ is a convex set with non-empty interior, we can find a vector $\xi \in \mathbb{R}^n$ such that $\langle N_i,\xi \rangle > 0$ for all $i \in I_0$. On the other hand, since $a_i = 0$ for all $i \in I \setminus I_0$, we obtain $\sum_{i \in I_0} a_i \, \langle N_i,\xi \rangle = \sum_{i \in I} a_i \, \langle N_i,\xi \rangle = 0$ and $\sum_{i \in I_0} a_i = \sum_{i \in I} a_i = 1$. This is a contradiction. This completes the proof of Lemma \ref{Lambda}. \\

\begin{lemma}
\label{transversality.1}
We can find a large constant $\Xi \in (2\Lambda,\infty)$ with the following property. Suppose that $1 \leq k \leq q$, and suppose that $x$ is a point in $\Omega$ with $-2\Xi^{-1} \leq u_k(x) \leq 0$. Moreover, suppose that $a_0,\hdots,a_{k-1}$ are nonnegative real numbers such that $a_i = 0$ for all $0 \leq i \leq k-1$ satisfying $u_i(x) \leq -2\Xi^{-1}$. Then 
\[\sum_{i=0}^{k-1} a_i \leq \Xi \, \Big | \Big ( \sum_{i=0}^{k-1} a_i \, N_i \Big ) \wedge N_k \Big |.\] 
\end{lemma} 

\textbf{Proof.} 
Let $\Lambda$ denote the constant in Lemma \ref{Lambda}. Let us fix a positive real number $\delta \in (0,\Lambda^{-1})$ so that the conclusion of Lemma \ref{approximate.angle.inequality.1} holds with $\varepsilon = \frac{1}{2} \, \Lambda^{-1}$. We define $\Xi = 2\delta^{-1} \in (2\Lambda,\infty)$. We claim that $\Xi$ has the desired property. To prove this, suppose that $1 \leq k \leq q$, and suppose that $x$ is a point in $\Omega$ with $-\delta \leq u_k(x) \leq 0$. Moreover, suppose that $a_0,\hdots,a_{k-1}$ are nonnegative real numbers such that $a_i = 0$ for all $0 \leq i \leq k-1$ satisfying $u_i(x) \leq -\delta$. We define 
\[b = -\Big \langle \sum_{i=0}^{k-1} a_i \, N_i,N_k \Big \rangle.\] 
Clearly, 
\[\Big | \Big ( \sum_{i=0}^{k-1} a_i \, N_i \Big ) \wedge N_k \Big | = \Big | \sum_{i=0}^{k-1} a_i \, N_i + b \, N_k \Big |.\] 
To estimate the term on the right hand side, we distinguish two cases:

\textit{Case 1:} Suppose first that $b \geq 0$. Lemma \ref{Lambda} implies 
\[\Big | \sum_{i=0}^{k-1} a_i \, N_i + b \, N_k \Big | \geq \Lambda^{-1} \, \Big ( \sum_{i=0}^{k-1} a_i + b \Big ) \geq \Lambda^{-1} \sum_{i=0}^{k-1} a_i.\] 

\textit{Case 2:} Suppose next that $b \leq 0$. Using Lemma \ref{Lambda}, we obtain 
\[\Big | \sum_{i=0}^{k-1} a_i \, N_i \Big | \geq \Lambda^{-1} \sum_{i=0}^{k-1} a_i.\] 
Moreover, it follows from Lemma \ref{approximate.angle.inequality.1} that $\langle N_i,N_k \rangle \leq \frac{1}{2} \, \Lambda^{-1}$ for all $0 \leq i \leq k-1$ satisfying $-\delta \leq u_i(x) \leq 0$. This implies 
\[b = -\sum_{i=0}^{k-1} a_i \, \langle N_i,N_k \rangle \geq -\frac{1}{2} \, \Lambda^{-1} \sum_{i=0}^{k-1} a_i.\] 
Consequently, 
\[\Big | \sum_{i=0}^{k-1} a_i \, N_i + b \, N_k \Big | \geq \Big | \sum_{i=0}^{k-1} a_i \, N_i \Big | + b \geq \frac{1}{2} \, \Lambda^{-1} \sum_{i=0}^{k-1} a_i.\] 
This completes the proof of Lemma \ref{transversality.1}. \\

Throughout this paper, we fix a smooth, even function $\eta: \mathbb{R} \to \mathbb{R}$ such that $\eta(t) = |t|$ for $|t| \geq \frac{1}{2}$ and $\eta''(t) \geq 0$ for all $t \in \mathbb{R}$. Note that $-1 \leq \eta'(t) \leq 0$ for $t \leq 0$ and $0 \leq \eta'(t) \leq 1$ for $t \geq 0$. Moreover, $|t| \leq \eta(t) \leq |t|+1$ for all $t \in \mathbb{R}$. 

In the following, $\gamma \in (0,\frac{1}{2})$ will denote a small parameter, and $\lambda_0 > 1$ will denote a large parameter. For abbreviation, we put $\lambda_k := \gamma^{-k} \lambda_0$ for $1 \leq k \leq q$.

\begin{definition}
We define a collection of smooth functions $\hat{u}_0,\hdots,\hat{u}_q$ on $\mathbb{R}^n$ so that $\hat{u}_0 = u_0$ and 
\[\hat{u}_k = \frac{1}{2} \, \big ( \hat{u}_{k-1} + u_k + \lambda_k^{-1} \, \eta(\lambda_k (\hat{u}_{k-1} - u_k)) \big )\] 
for $1 \leq k \leq q$. Moreover, we define $\hat{\Omega} := \{\hat{u}_q \leq 0\}$ and $\hat{\Sigma} := \{\hat{u}_q = 0\}$.
\end{definition}

\begin{lemma}
\label{inequality.for.hat.u_k}
We have 
\[\max \{\hat{u}_{k-1},u_k\} \leq \hat{u}_k \leq \max \{\hat{u}_{k-1},u_k\} + \lambda_k^{-1}\] 
for $1 \leq k \leq q$. In particular, $\hat{u}_0 \leq \hat{u}_1 \leq \hdots \leq \hat{u}_q$. 
\end{lemma} 

\textbf{Proof.} This follows from the fact that $|t| \leq \eta(t) \leq |t|+1$ for all $t \in \mathbb{R}$. \\

\begin{lemma}
\label{inclusion}
We have 
\[\max \{u_0,\hdots,u_q\} \leq \hat{u}_q \leq \max \{u_0,\hdots,u_q\} + \sum_{k=1}^q \lambda_k^{-1}.\] 
In particular, $\bigcap_{m=0}^q \{u_m \leq -\lambda_0^{-1}\} \subset \hat{\Omega} \subset \Omega$. 
\end{lemma}

\textbf{Proof.} This follows immediately from Lemma \ref{inequality.for.hat.u_k}. \\

\begin{lemma}
\label{C2.bound.for.hat.u}
Let $0 \leq k \leq q$. Then $|d\hat{u}_k| \leq C$ at each point in $\mathbb{R}^n$. Moreover, $\hat{u}_k$ is a convex function and the Hessian of $\hat{u}_k$ is bounded by $C\lambda_k$ at each point in $\mathbb{R}^n$. Here, $C$ is independent of $\gamma$ and $\lambda_0$.
\end{lemma}

\textbf{Proof.} 
The proof is by induction on $k$. The assertion is obvious for $k=0$. Suppose next that $1 \leq k \leq q$ and that the assertion is true for $\hat{u}_{k-1}$. The differential of $u_k$ is given by 
\[d\hat{u}_k = \frac{1}{2} \, (1+\eta'(\lambda_k(\hat{u}_{k-1}-u_k))) \, d\hat{u}_{k-1} + \frac{1}{2} \, (1-\eta'(\lambda_k(\hat{u}_{k-1}-u_k))) \, du_k.\] 
Moreover, 
\begin{align*} 
D^2 \hat{u}_k 
&= \frac{1}{2} \, (1+\eta'(\lambda_k(\hat{u}_{k-1}-u_k))) \, D^2 \hat{u}_{k-1} \\ 
&+ \frac{1}{2} \, \lambda_k \, \eta''(\lambda_k(\hat{u}_{k-1}-u_k)) \, (d\hat{u}_{k-1} - du_k) \otimes (d\hat{u}_{k-1} - du_k), 
\end{align*}
where $D^2 \hat{u}_k$ denotes the Hessian of $\hat{u}_k$ with respect to the Euclidean metric and $D^2 \hat{u}_{k-1}$ denotes the Hessian of $\hat{u}_{k-1}$ with respect to the Euclidean metric. From this, we deduce that the assertion is true for $\hat{u}_k$. This completes the proof of Lemma \ref{C2.bound.for.hat.u}. \\

\begin{lemma}
\label{gradient.of.hat.u_k}
Let $0 \leq k \leq q$ and let $x$ be a point in $\Omega$. Then we can find nonnegative real numbers $a_0,\hdots,a_k$ such that $\sum_{i=0}^k a_i = 1$ and $d\hat{u}_k = \sum_{i=0}^k a_i \, du_i$ at the point $x$. Moreover, $a_i=0$ for all $0 \leq i \leq k$ satisfying $u_i(x) < \hat{u}_k(x) - 2k\lambda_i^{-1}$.
\end{lemma} 

\textbf{Proof.} 
We argue by induction on $k$. The assertion is clearly true for $k=0$. We next assume that $1 \leq k \leq q$, and the assertion is true for $k-1$. To prove the assertion for $k$, we fix an arbitrary point $x \in \Omega$.

\textit{Case 1:} Suppose that $\hat{u}_{k-1}(x)-u_k(x) > \lambda_k^{-1}$. In this case, $\hat{u}_k = \hat{u}_{k-1}$ in an open neighborhood of $x$. In view of the induction hypothesis, we can find nonnegative real numbers $b_0,\hdots,b_{k-1}$ such that $\sum_{i=0}^{k-1} b_i = 1$ and $d\hat{u}_{k-1} = \sum_{i=0}^{k-1} b_i \, du_i$ at the point $x$. Moreover, $b_i=0$ for all $0 \leq i \leq k-1$ satisfying $u_i(x) < \hat{u}_{k-1}(x) - 2(k-1) \lambda_i^{-1}$. We define nonnegative real numbers $a_0,\hdots,a_k$ by $a_i = b_i$ for $0 \leq i \leq k-1$ and $a_k = 0$. Then $\sum_{i=0}^k a_i = \sum_{i=0}^{k-1} b_i = 1$ and $d\hat{u}_k = d\hat{u}_{k-1} = \sum_{i=0}^{k-1} b_i \, du_i = \sum_{i=0}^k a_i \, du_i$ at the point $x$. Moreover, $a_i=0$ for all $0 \leq i \leq k$ satisfying $u_i(x) < \hat{u}_k(x) - 2k \lambda_i^{-1}$. 

\textit{Case 2:} Suppose that $\hat{u}_{k-1}(x)-u_k(x) < -\lambda_k^{-1}$. In this case, $\hat{u}_k = u_k$ in an open neighborhood of $x$. We define nonnegative real numbers $a_0,\hdots,a_k$ by $a_i=0$ for $0 \leq i \leq k-1$ and $a_k=1$. Then $\sum_{i=0}^k a_i = 1$ and $d\hat{u}_k = du_k = \sum_{i=0}^k a_i \, du_i$. Moreover, $a_i=0$ for all $0 \leq i \leq k$ satisfying $u_i(x) < \hat{u}_k(x) - 2k \lambda_i^{-1}$. 

\textit{Case 3:} Suppose that $|\hat{u}_{k-1}(x)-u_k(x)| \leq \lambda_k^{-1}$. In this case, Lemma \ref{inequality.for.hat.u_k} implies 
\[\hat{u}_k(x) \leq \max \{\hat{u}_{k-1}(x),u_k(x)\} + \lambda_k^{-1} \leq \min \{\hat{u}_{k-1}(x),u_k(x)\} + 2\lambda_k^{-1}.\] 
In view of the induction hypothesis, we can find nonnegative real numbers $b_0,\hdots,b_{k-1}$ such that $\sum_{i=0}^{k-1} b_i = 1$ and $d\hat{u}_{k-1} = \sum_{i=0}^{k-1} b_i \, du_i$ at the point $x$. Moreover, $b_i=0$ for all $0 \leq i \leq k-1$ satisfying $u_i(x) < \hat{u}_{k-1}(x) - 2(k-1) \lambda_i^{-1}$. We define nonnegative real numbers $a_0,\hdots,a_k$ by 
\[a_i = \frac{1}{2} \, (1+\eta'(\lambda_k(\hat{u}_{k-1}-u_k))) \, b_i\] 
for $0 \leq i \leq k-1$, and 
\[a_k = \frac{1}{2} \, (1-\eta'(\lambda_k(\hat{u}_{k-1}-u_k))).\] 
Then 
\begin{align*} 
\sum_{i=0}^k a_i 
&= \frac{1}{2} \, (1+\eta'(\lambda_k(\hat{u}_{k-1}-u_k))) \, \sum_{i=0}^{k-1} b_i + \frac{1}{2} \, (1-\eta'(\lambda_k(\hat{u}_{k-1}-u_k))) \\ 
&= \frac{1}{2} \, (1+\eta'(\lambda_k(\hat{u}_{k-1}-u_k))) + \frac{1}{2} \, (1-\eta'(\lambda_k(\hat{u}_{k-1}-u_k))) \\ 
&= 1 
\end{align*} 
and 
\begin{align*} 
d\hat{u}_k 
&= \frac{1}{2} \, (1+\eta'(\lambda_k(\hat{u}_{k-1}-u_k))) \, d\hat{u}_{k-1} + \frac{1}{2} \, (1-\eta'(\lambda_k(\hat{u}_{k-1}-u_k))) \, du_k \\ 
&= \frac{1}{2} \, (1+\eta'(\lambda_k(\hat{u}_{k-1}-u_k))) \, \sum_{i=0}^{k-1} b_i \, du_i + \frac{1}{2} \, (1-\eta'(\lambda_k(\hat{u}_{k-1}-u_k))) \, du_k \\ 
&= \sum_{i=0}^k a_i \, du_i 
\end{align*}
at the point $x$. Finally, $a_i=0$ for all $0 \leq i \leq k$ satisfying $u_i(x) < \hat{u}_k(x) - 2k \lambda_i^{-1}$. This completes the proof of Lemma \ref{gradient.of.hat.u_k}. \\

\begin{lemma}
\label{lower.bound.for.norm.of.gradient.of.hat.u_k}
If $\lambda_0$ is sufficiently large, then the following holds. Let $0 \leq k \leq q$. If $x$ is a point in $\Omega$ satisfying $-\Lambda^{-1} \leq \hat{u}_k(x) \leq 0$, then $|d\hat{u}_k| \geq \Lambda^{-1}$ at the point $x$. In particular, $\hat{\Sigma}$ is a smooth hypersurface.
\end{lemma}

\textbf{Proof.} 
This follows by combining Lemma \ref{Lambda} and Lemma \ref{gradient.of.hat.u_k}. \\

\begin{lemma}
\label{transversality.2}
If $\lambda_0$ is sufficiently large, then the following holds. Let $0 \leq j < k \leq q$. If $x$ is a point in $\Omega$ satisfying $-\Xi^{-1} \leq \hat{u}_j(x) \leq 0$ and $-\Xi^{-1} \leq u_k(x) \leq 0$, then $|d\hat{u}_j \wedge du_k| \geq \Xi^{-1}$ at the point $x$. 
\end{lemma}

\textbf{Proof.} 
Consider a point $x \in \Omega$ satisfying $-\Xi^{-1} \leq \hat{u}_j(x) \leq 0$ and $-\Xi^{-1} \leq u_k(x) \leq 0$. By Lemma \ref{gradient.of.hat.u_k}, we can find nonnegative real numbers $a_0,\hdots,a_j$ such that $\sum_{i=0}^j a_i = 1$ and $d\hat{u}_j = \sum_{i=0}^j a_i \, du_i$ at the point $x$. Moreover, $a_i=0$ for all $i \in \{0,1,\hdots,j\}$ satisfying $u_i(x) < \hat{u}_j(x) - 2j \lambda_i^{-1}$. Since $-\Xi^{-1} \leq \hat{u}_j(x) \leq 0$, it follows that $a_i=0$ for all $i \in \{0,1,\hdots,j\}$ satisfying $u_i(x) \leq -2\Xi^{-1}$. Using Assumption \ref{gradient.of.u_k} and Lemma \ref{transversality.1}, we obtain 
\[1 = \sum_{i=0}^j a_i \leq \Xi \, \Big | \Big ( \sum_{i=0}^j a_i \, du_i \Big ) \wedge du_k \Big | = \Xi \, |d\hat{u}_j \wedge du_k|\] 
at the point $x$. This completes the proof of Lemma \ref{transversality.2}. \\

\begin{lemma}
\label{area.1}
Let $1 \leq k \leq q$. Then 
\begin{align*} 
&|\Omega \cap \{\hat{u}_k=0\} \cap \{-2\lambda_k^{-1} \leq \hat{u}_{k-1} \leq 0\} \cap \{-2\lambda_k^{-1} \leq u_k \leq 0\} \cap B_r(p)| \\ 
&\leq C \lambda_k^{-1} r^{n-2} 
\end{align*} 
for all $0 < r \leq 1$. Here, the expression on the left hand side represents the $(n-1)$-dimensional measure. Moreover, $B_r(p)$ denotes a Euclidean ball of radius $r$, and $C$ is independent of $\gamma$ and $\lambda_0$.
\end{lemma} 

\textbf{Proof.} 
Lemma \ref{transversality.2} implies that $|d\hat{u}_{k-1} \wedge du_k| \geq \Xi^{-1}$ at each point in $\Omega \cap \{-2\lambda_k^{-1} \leq \hat{u}_{k-1} \leq 0\} \cap \{-2\lambda_k^{-1} \leq u_k \leq 0\}$. Let us consider two arbitrary real numbers $a,b \in [-2\lambda_k^{-1},0]$. We apply Proposition \ref{area.estimate} with $m = n-1$ and with $f$ defined as the restriction of the function $\hat{u}_{k-1}-a$ to the hyperplane $\{u_k=b\}$. This gives 
\[|\{d\hat{u}_{k-1} \wedge du_k \neq 0\} \cap \{\hat{u}_{k-1}=a\} \cap \{u_k=b\} \cap B_r(p)| \leq C(n) \, r^{n-2}\] 
for all $0 < r \leq 1$, where the expression on the left hand side represents the $(n-2)$-dimensional measure. Since $d\hat{u}_{k-1} \wedge du_k \neq 0$ at each point in $\Omega \cap \{-2\lambda_k^{-1} \leq \hat{u}_{k-1} \leq 0\} \cap \{-2\lambda_k^{-1} \leq u_k \leq 0\}$, it follows that 
\begin{equation} 
\label{area.bound}
|\Omega \cap \{\hat{u}_{k-1}=a\} \cap \{u_k=b\} \cap B_r(p)| \leq C(n) \,  r^{n-2} 
\end{equation}
for all $a,b \in [-2\lambda_k^{-1},0]$ and all $0 < r \leq 1$, where the expression on the left hand side represents the $(n-2)$-dimensional measure. 

For abbreviation, let 
\[S_k = \Omega \cap \{\hat{u}_k=0\} \cap \{-2\lambda_k^{-1} \leq \hat{u}_{k-1} \leq 0\} \cap \{-2\lambda_k^{-1} \leq u_k \leq 0\}.\] 
Moreover, we put $S_k^+ = S_k \cap \{\hat{u}_{k-1}-u_k \geq 0\}$ and $S_k^- = S_k \cap \{\hat{u}_{k-1}-u_k \leq 0\}$. 

For each $b \in [-2\lambda_k^{-1},0]$, we can find a real number $a \in [b,0]$ (depending on $b$) such that $S_k^+ \cap \{u_k=b\} \subset \{\hat{u}_{k-1}=a\}$. Using (\ref{area.bound}), we obtain 
\begin{equation} 
\label{intersection.of.S_k.with.level.sets.of.u_k}
|S_k^+ \cap \{u_k=b\} \cap B_r(p)| \leq C(n) \, r^{n-2} 
\end{equation}
for each $b \in [-2\lambda_k^{-1},0]$ and all $0 < r \leq 1$, where the expression on the left hand side represents the $(n-2)$-dimensional measure. Similarly, for each $a \in [-2\lambda_k^{-1},0]$, we can find a real number $b \in [a,0]$ (depending on $a$) such that $S_k^- \cap \{\hat{u}_{k-1}=a\} \subset \{u_k=b\}$. Using (\ref{area.bound}), we obtain 
\begin{equation} 
\label{intersection.of.S_k.with.level.sets.of.hat.u_k-1}
|S_k^- \cap \{\hat{u}_{k-1}=a\} \cap B_r(p)| \leq C(n) \, r^{n-2} 
\end{equation}
for each $a \in [-2\lambda_k^{-1},0]$ and all $0 < r \leq 1$, where the expression on the left hand side represents the $(n-2)$-dimensional measure. In the next step, we integrate the inequality (\ref{intersection.of.S_k.with.level.sets.of.u_k}) over $b \in [-2\lambda_k^{-1},0]$, and we integrate the inequality (\ref{intersection.of.S_k.with.level.sets.of.hat.u_k-1}) over $a \in [-2\lambda_k^{-1},0]$. Using the co-area formula, we obtain 
\[\int_{S_k^+ \cap B_r(p)} \frac{|du_k \wedge d\hat{u}_k|}{|d\hat{u}_k|} \leq C(n) \, \lambda_k^{-1} r^{n-2}\] 
and 
\[\int_{S_k^- \cap B_r(p)} \frac{|d\hat{u}_{k-1} \wedge d\hat{u}_k|}{|d\hat{u}_k|} \leq C(n) \, \lambda_k^{-1} r^{n-2}\] 
for all $0 < r \leq 1$. It follows from Lemma \ref{C2.bound.for.hat.u} that $|d\hat{u}_k| \leq C$ at each point in $\Omega$. Using the identity 
\[d\hat{u}_k = \frac{1}{2} \, (1+\eta'(\lambda_k(\hat{u}_{k-1}-u_k))) \, d\hat{u}_{k-1} + \frac{1}{2} \, (1-\eta'(\lambda_k(\hat{u}_{k-1}-u_k))) \, du_k,\] 
we obtain 
\begin{align*} 
|du_k \wedge d\hat{u}_k| 
&= \frac{1}{2} \, (1+\eta'(\lambda_k(\hat{u}_{k-1}-u_k))) \, |du_k \wedge d\hat{u}_{k-1}| \\ 
&\geq \frac{1}{2} \, |du_k \wedge d\hat{u}_{k-1}| \geq \frac{1}{2} \, \Xi^{-1} 
\end{align*} 
at each point in $S_k^+$ and
\begin{align*} 
|d\hat{u}_{k-1} \wedge d\hat{u}_k| 
&= \frac{1}{2} \, (1-\eta'(\lambda_k(\hat{u}_{k-1}-u_k))) \, |d\hat{u}_{k-1} \wedge du_k| \\ 
&\geq \frac{1}{2} \, |d\hat{u}_{k-1} \wedge du_k| \geq \frac{1}{2} \, \Xi^{-1} 
\end{align*}
at each point in $S_k^-$. Putting these facts together, we conclude that 
\[|S_k^+ \cap B_r(p)| \leq C \lambda_k^{-1} r^{n-2}\] 
and 
\[|S_k^- \cap B_r(p)| \leq C \lambda_k^{-1} r^{n-2}\] 
for all $0 < r \leq 1$. Thus, $|S_k \cap B_r(p)| \leq C \lambda_k^{-1} r^{n-2}$ for all $0 < r \leq 1$. This completes the proof of Lemma \ref{area.1}. \\

\begin{lemma}
\label{area.2}
Suppose that $i,j,k \in \{0,1,\hdots,q\}$ are pairwise distinct. If $\lambda_0$ is sufficiently large, then 
\begin{align*} 
&|\Omega \cap \{\hat{u}_k=0\} \cap \{-2\lambda_k^{-1} \leq \hat{u}_{k-1} \leq 0\} \cap \{-2\lambda_k^{-1} \leq u_k \leq 0\} \\ 
&\cap \{-4\lambda_0^{-1} \leq u_j \leq 0\} \cap \{-6\lambda_0^{-1} \leq u_i \leq 0\} \cap B_r(p)| \leq C \lambda_k^{-1} \lambda_0^{-1} r^{n-3} 
\end{align*}
for all $\lambda_0^{-1} \leq r \leq 1$. Here, the expression on the left hand side represents the $(n-1)$-dimensional measure. Moreover, $B_r(p)$ denotes a Euclidean ball of radius $r$, and $C$ is independent of $\gamma$ and $\lambda_0$. 
\end{lemma}

\textbf{Proof.} As above, we define 
\[S_k = \Omega \cap \{\hat{u}_k=0\} \cap \{-2\lambda_k^{-1} \leq \hat{u}_{k-1} \leq 0\} \cap \{-2\lambda_k^{-1} \leq u_k \leq 0\}.\] 
We distinguish two cases: 

\textit{Case 1:} Suppose that $\Omega \cap \{u_i=0\} \cap \{u_j=0\} \cap \{u_k=0\} = \emptyset$. Then $\Omega \cap \{-2\lambda_0^{-1} \leq u_k \leq 0\} \cap \{-4\lambda_0^{-1} \leq u_j \leq 0\} \cap \{-6\lambda_0^{-1} \leq u_i \leq 0\} = \emptyset$ if $\lambda_0$ is sufficiently large. Consequently, $S_k \cap \{-4\lambda_0^{-1} \leq u_j \leq 0\} \cap \{-6\lambda_0^{-1} \leq u_i \leq 0\} = \emptyset$ if $\lambda_0$ is sufficiently large. Thus, the assertion is trivially true in this case.

\textit{Case 2:} Suppose that $\Omega \cap \{u_i=0\} \cap \{u_j=0\} \cap \{u_k=0\} \neq \emptyset$. It follows from Assumption \ref{no.redundant.inequalities} that the hyperplanes $\{u_i=0\}$, $\{u_j=0\}$, $\{u_k=0\}$ must intersect transversally. Suppose now that $\lambda_0^{-1} \leq r \leq 1$. We can cover the set $\{-2\lambda_0^{-1} \leq u_k \leq 0\} \cap \{-4\lambda_0^{-1} \leq u_j \leq 0\} \cap \{-6\lambda_0^{-1} \leq u_i \leq 0\} \cap B_r(p)$ by $C \, (\lambda_0 r)^{n-3}$ balls of radius $\lambda_0^{-1}$. By Lemma \ref{area.1}, the intersection of $S_k$ with each such ball has area at most $C \lambda_k^{-1} \lambda_0^{2-n}$. Therefore, the set $S_k \cap \{-4\lambda_0^{-1} \leq u_j \leq 0\} \cap \{-6\lambda_0^{-1} \leq u_i \leq 0\} \cap B_r(p)$ has area at most $C \lambda_k^{-1} \lambda_0^{-1} r^{n-3}$. This completes the proof of Lemma \ref{area.2}. \\

\begin{definition}
We define 
\[F_0 = \hat{\Sigma} \cap \bigcap_{m=1}^q \{\hat{u}_{m-1}-u_m > \lambda_m^{-1}\}.\] 
For $1 \leq k \leq q$, we define 
\[F_k = \hat{\Sigma} \cap \bigcap_{m=k+1}^q \{\hat{u}_{m-1}-u_m > \lambda_m^{-1}\} \cap \{\hat{u}_{k-1}-u_k < -\lambda_k^{-1}\}.\]  
\end{definition}

\begin{definition}
For $1 \leq k \leq q$, we define 
\begin{align*} 
E_{0,k} 
&= \hat{\Sigma} \cap \bigcap_{m=k+1}^q \{\hat{u}_{m-1}-u_m > \lambda_m^{-1}\} \cap \bigcap_{m=1}^{k-1} \{\hat{u}_{m-1}-u_m > \lambda_m^{-1}\} \\ 
&\cap \{-2\lambda_k^{-1} \leq \hat{u}_{k-1} \leq 0\} \cap \{-2\lambda_k^{-1} \leq u_k \leq 0\} \cap \{-2\lambda_k^{-1} \leq u_0 \leq 0\}. 
\end{align*}
For $1 \leq j < k \leq q$, we define 
\begin{align*} 
E_{j,k} 
&= \hat{\Sigma} \cap \bigcap_{m=k+1}^q \{\hat{u}_{m-1}-u_m > \lambda_m^{-1}\} \cap \bigcap_{m=j+1}^{k-1} \{\hat{u}_{m-1}-u_m > \lambda_m^{-1}\} \\ 
&\cap \{-2\lambda_k^{-1} \leq \hat{u}_{k-1} \leq 0\} \cap \{-2\lambda_k^{-1} \leq u_k \leq 0\} \\ 
&\cap \{-2\lambda_k^{-1} \leq u_j \leq 0\} \cap \{\hat{u}_{j-1}-u_j < -\lambda_j^{-1}\}.
\end{align*} 
\end{definition}

\begin{definition} 
For $0 \leq i < j < k \leq q$, we define 
\begin{align*} 
G_{i,j,k} 
&= \hat{\Sigma} \cap \bigcap_{m=k+1}^q \{\hat{u}_{m-1}-u_m > \lambda_m^{-1}\} \cap \bigcap_{m=j+1}^{k-1} \{\hat{u}_{m-1}-u_m > \lambda_m^{-1}\} \\ 
&\cap \{-2\lambda_k^{-1} \leq \hat{u}_{k-1} \leq 0\} \cap \{-2\lambda_k^{-1} \leq u_k \leq 0\} \\ 
&\cap \{-4\lambda_j^{-1} \leq u_j \leq 0\} \cap \{-6\lambda_i^{-1} \leq u_i \leq 0\}.
\end{align*}
\end{definition}

\begin{proposition}
\label{decomposition.into.subsets} 
We have
\[\hat{\Sigma} = \bigcup_{0 \leq k \leq q} F_k \cup \bigcup_{0 \leq j < k \leq q} E_{j,k} \cup \bigcup_{0 \leq i<j<k \leq q} G_{i,j,k}.\]
\end{proposition} 

\textbf{Proof.} Let us consider an arbitrary point $x \in \hat{\Sigma}$. We claim that 
\[x \in \bigcup_{0 \leq k \leq q} F_k \cup \bigcup_{0 \leq j < k \leq q} E_{j,k} \cup \bigcup_{0 \leq i<j<k \leq q} G_{i,j,k}.\] 
If $\hat{u}_{m-1}(x) - u_m(x) > \lambda_m^{-1}$ for all $m \in \{1,\hdots,q\}$, then $x \in F_0$. Otherwise, we define an integer $k \in \{1,\hdots,q\}$ by 
\[k = \max \big \{ m \in \{1,\hdots,q\}: \hat{u}_{m-1}(x) - u_m(x) \leq \lambda_m^{-1} \big \}.\] 
Clearly, $\hat{u}_{k-1}(x) - u_k(x) \leq \lambda_k^{-1}$ and $\hat{u}_{m-1}(x) - u_m(x) > \lambda_m^{-1}$ for all $m \in \{k+1,\hdots,q\}$. From this, we deduce that $\hat{u}_k(x) = \hat{u}_q(x) = 0$. If $\hat{u}_{k-1}(x)-u_k(x) < -\lambda_k^{-1}$, then $x \in F_k$.

It remains to consider the case when $|\hat{u}_{k-1}(x)-u_k(x)| \leq \lambda_k^{-1}$. Since $\hat{u}_k(x) = 0$, Lemma \ref{inequality.for.hat.u_k} implies that $\max \{\hat{u}_{k-1}(x),u_k(x)\} \in [-\lambda_k^{-1},0]$. Therefore, $\hat{u}_{k-1}(x) \in [-2\lambda_k^{-1},0]$ and $u_k(x) \in [-2\lambda_k^{-1},0]$. If $\hat{u}_{m-1}(x) - u_m(x) > \lambda_m^{-1}$ for all $m \in \{1,\hdots,k-1\}$, then $u_0(x) = \hat{u}_0(x) = \hat{u}_{k-1}(x) \in [-2\lambda_k^{-1},0]$, and consequently $x \in E_{0,k}$. Otherwise, we define an integer $j \in \{1,\hdots,k-1\}$ by 
\[j = \max \big \{ m \in \{1,\hdots,k-1\}: \hat{u}_{m-1}(x) - u_m(x) \leq \lambda_m^{-1} \big \}.\] 
Clearly, $\hat{u}_{j-1}(x) - u_j(x) \leq \lambda_j^{-1}$ and $\hat{u}_{m-1}(x) - u_m(x) > \lambda_m^{-1}$ for all $m \in \{j+1,\hdots,k-1\}$. From this, we deduce that $\hat{u}_j(x) = \hat{u}_{k-1}(x) \in [-2\lambda_k^{-1},0]$. If $\hat{u}_{j-1}(x) - u_j(x) < -\lambda_j^{-1}$, then $u_j(x) = \hat{u}_j(x) = \hat{u}_{k-1}(x) \in [-2\lambda_k^{-1},0]$, and consequently $x \in E_{j,k}$. 

It remains to consider the case when $|\hat{u}_{j-1}(x)-u_j(x)| \leq \lambda_j^{-1}$. Since $\hat{u}_j(x) \in [-2\lambda_k^{-1},0]$, Lemma \ref{inequality.for.hat.u_k} implies that $\max \{\hat{u}_{j-1}(x),u_j(x)\} \in [-3\lambda_j^{-1},0]$. Therefore, $\hat{u}_{j-1}(x) \in [-4\lambda_j^{-1},0]$ and $u_j(x) \in [-4\lambda_j^{-1},0]$. If $\hat{u}_{m-1}(x) - u_m(x) > \lambda_m^{-1}$ for all $m \in \{1,\hdots,j-1\}$, then $u_0(x) = \hat{u}_0(x) = \hat{u}_{j-1}(x) \in [-4\lambda_j^{-1},0]$, and consequently $x \in G_{0,j,k}$. Otherwise, we define an integer $i \in \{1,\hdots,j-1\}$ by 
\[i = \max \big \{ m \in \{1,\hdots,j-1\}: \hat{u}_{m-1}(x) - u_m(x) \leq \lambda_m^{-1} \big \}.\] 
Clearly, $\hat{u}_{i-1}(x) - u_i(x) \leq \lambda_i^{-1}$ and $\hat{u}_{m-1}(x) - u_m(x) > \lambda_m^{-1}$ for all $m \in \{i+1,\hdots,j-1\}$. From this, we deduce that $\hat{u}_i(x) = \hat{u}_{j-1}(x) \in [-4\lambda_j^{-1},0]$. Using Lemma \ref{inequality.for.hat.u_k}, we obtain $\max \{\hat{u}_{i-1}(x),u_i(x)\} \in [-5\lambda_i^{-1},0]$. Since $\hat{u}_{i-1}(x) - u_i(x) \leq \lambda_i^{-1}$, we conclude that $u_i(x) \in [-6\lambda_i^{-1},0]$, and consequently $x \in G_{i,j,k}$. This completes the proof of Proposition \ref{decomposition.into.subsets}. \\

\section{The map $\hat{N}: \hat{\Sigma} \to S^{n-1}$}

\label{construction.hat.N}

Suppose that $n \geq 3$ is an integer, and $\Omega = \bigcap_{m=0}^q \{u_m \leq 0\}$ is a compact, convex polytope in $\mathbb{R}^n$ with non-empty interior. Throughout this section, we assume that Assumptions \ref{no.redundant.inequalities}, \ref{gradient.of.u_k}, and \ref{angles.bounded.by.pi/2} are satisfied. Moreover, we assume that $g$ is an arbitrary Riemannian metric on $\mathbb{R}^n$.

\begin{definition}
Let $0 \leq k \leq q$. At each point $x \in \mathbb{R}^n$, we define $\nu_k = \frac{\nabla u_k}{|\nabla u_k|}$, where $\nabla u_k$ and $|\nabla u_k|$ are computed with respect to the metric $g$. At each point $x \in \mathbb{R}^n \setminus \{d\hat{u}_k = 0\}$, we define $\hat{\nu}_k = \frac{\nabla \hat{u}_k}{|\nabla \hat{u}_k|}$, where $\nabla \hat{u}_k$ and $|\nabla \hat{u}_k|$ are computed with respect to the metric $g$. Note that $\Omega \cap \{-2^{-k} \, \Xi^{-1} \leq \hat{u}_k \leq 0\} \subset \mathbb{R}^n \setminus \{d\hat{u}_k = 0\}$ by Lemma \ref{lower.bound.for.norm.of.gradient.of.hat.u_k}. Finally, we denote by $\hat{\nu}$ the restriction of the vector field $\hat{\nu}_q$ to the hypersurface $\hat{\Sigma} = \{\hat{u}_q=0\}$.
\end{definition}

\begin{proposition}
\label{recursive.relation.for.hat.nu}
If $\lambda_0$ is sufficiently large, then the following statements hold: \\ 
(i) Let $1 \leq k \leq q$. Suppose that $x$ is a point in $\Omega$ satisfying $-2^{-k} \, \Xi^{-1} \leq \hat{u}_k(x) \leq 0$ and $\hat{u}_{k-1}(x)-u_k(x) > \lambda_k^{-1}$. Then $-2^{-k+1} \, \Xi^{-1} \leq \hat{u}_{k-1}(x) \leq 0$. Moreover, $\hat{\nu}_k = \hat{\nu}_{k-1}$ at the point $x$. \\
(ii) Let $1 \leq k \leq q$. Suppose that $x$ is a point in $\Omega$ satisfying $-2^{-k} \, \Xi^{-1} \leq \hat{u}_k(x) \leq 0$ and $\hat{u}_{k-1}(x)-u_k(x) < -\lambda_k^{-1}$. Then $\hat{\nu}_k = \nu_k$ at the point $x$. \\
(iii) Let $1 \leq k \leq q$. Suppose that $x$ is a point in $\Omega$ satisfying $-2^{-k} \, \Xi^{-1} \leq \hat{u}_k(x) \leq 0$ and $|\hat{u}_{k-1}(x)-u_k(x)| \leq \lambda_k^{-1}$. Then $-2^{-k+1} \, \Xi^{-1} \leq \hat{u}_{k-1}(x) \leq 0$. We can find a real number $\alpha \in (0,\frac{\pi}{2})$ such that $\cos(2\alpha) = \langle \hat{\nu}_{k-1},\nu_k \rangle$ at the point $x$, where the inner product is computed with respect to the metric $g$. Let us define a real number $\varphi \in [-\alpha,\alpha]$ by 
\begin{align*} 
&\frac{\tan(\varphi)}{\tan(\alpha)} \\ 
&= \frac{(1+\eta'(\lambda_k(\hat{u}_{k-1}-u_k))) \, |\nabla \hat{u}_{k-1}| - (1-\eta'(\lambda_k(\hat{u}_{k-1}-u_k))) \, |\nabla u_k|}{(1+\eta'(\lambda_k(\hat{u}_{k-1}-u_k))) \, |\nabla \hat{u}_{k-1}| + (1-\eta'(\lambda_k(\hat{u}_{k-1}-u_k))) \, |\nabla u_k|}, 
\end{align*} 
where $|\nabla \hat{u}_{k-1}|$ and $|\nabla u_k|$ are computed with respect to the metric $g$. Then 
\[\hat{\nu}_k = \frac{\sin(\alpha+\varphi) \, \hat{\nu}_{k-1} + \sin(\alpha-\varphi) \, \nu_k}{\sin(2\alpha)}\] 
at the point $x$.
\end{proposition}

\textbf{Proof.} 
Statement (i) and statement (ii) follow directly from the definition. 

To prove statement (iii), we consider an integer $1 \leq k \leq q$ and a point $x \in \Omega$ satisfying $-2^{-k} \, \Xi^{-1} \leq \hat{u}_k(x) \leq 0$ and $|\hat{u}_{k-1}(x)-u_k(x)| \leq \lambda_k^{-1}$. Using Lemma \ref{inequality.for.hat.u_k}, we obtain 
\[\max \{\hat{u}_{k-1}(x),u_k(x)\} \leq \hat{u}_k(x) \leq 0\] 
and 
\begin{align*} 
\min \{\hat{u}_{k-1}(x),u_k(x)\} 
&\geq \max \{\hat{u}_{k-1}(x),u_k(x)\} - \lambda_k^{-1} \\ 
&\geq \hat{u}_k(x) - 2\lambda_k^{-1} \\ 
&\geq -2^{-k+1} \, \Xi^{-1}. 
\end{align*}
Therefore, $-2^{-k+1} \, \Xi^{-1} \leq \hat{u}_{k-1}(x) \leq 0$ and $-2^{-k+1} \, \Xi^{-1} \leq u_k(x) \leq 0$. It follows from Lemma \ref{transversality.2} that $d\hat{u}_{k-1} \wedge du_k \neq 0$ at the point $x$. Consequently, we can find a real number $\alpha \in (0,\frac{\pi}{2})$ such that $\cos(2\alpha) = \langle \hat{\nu}_{k-1},\nu_k \rangle$ at the point $x$. We define a real number $\varphi \in [-\alpha,\alpha]$ so that 
\begin{align*} 
&\frac{\tan(\varphi)}{\tan(\alpha)} \\ 
&= \frac{(1+\eta'(\lambda_k(\hat{u}_{k-1}-u_k))) \, |\nabla \hat{u}_{k-1}| - (1-\eta'(\lambda_k(\hat{u}_{k-1}-u_k))) \, |\nabla u_k|}{(1+\eta'(\lambda_k(\hat{u}_{k-1}-u_k))) \, |\nabla \hat{u}_{k-1}| + (1-\eta'(\lambda_k(\hat{u}_{k-1}-u_k))) \, |\nabla u_k|} 
\end{align*} 
at the point $x$. Since $d\hat{u}_{k-1} \wedge du_k \neq 0$ at the point $x$, the denominator is strictly positive; this ensures that $\varphi$ is well-defined. 

We compute 
\begin{align*} 
&\frac{\sin(\alpha+\varphi) \, \hat{\nu}_{k-1} + \sin(\alpha-\varphi) \, \nu_k}{\sin(2\alpha)} \\ 
&= \frac{1}{2} \, \Big ( \frac{\cos(\varphi)}{\cos(\alpha)} + \frac{\sin(\varphi)}{\sin(\alpha)} \Big ) \, \hat{\nu}_{k-1} + \frac{1}{2} \, \Big ( \frac{\cos(\varphi)}{\cos(\alpha)} - \frac{\sin(\varphi)}{\sin(\alpha)} \Big ) \, \nu_k \\ 
&= \frac{1}{2} \, \frac{\cos(\varphi)}{\cos(\alpha)} \, \Big ( 1 + \frac{\tan(\varphi)}{\tan(\alpha)} \Big ) \, \frac{\nabla \hat{u}_{k-1}}{|\nabla \hat{u}_{k-1}|} + \frac{1}{2} \, \frac{\cos(\varphi)}{\cos(\alpha)} \, \Big ( 1 - \frac{\tan(\varphi)}{\tan(\alpha)} \Big ) \, \frac{\nabla u_k}{|\nabla u_k|} \\ 
&= \frac{\cos(\varphi)}{\cos(\alpha)} \, \frac{(1+\eta'(\lambda_k(\hat{u}_{k-1}-u_k))) \, \nabla \hat{u}_{k-1} + (1-\eta'(\lambda_k(\hat{u}_{k-1}-u_k))) \, \nabla u_k}{(1+\eta'(\lambda_k(\hat{u}_{k-1}-u_k))) \, |\nabla \hat{u}_{k-1}| + (1-\eta'(\lambda_k(\hat{u}_{k-1}-u_k))) \, |\nabla u_k|} \\ 
&= \frac{\cos(\varphi)}{\cos(\alpha)} \, \frac{2 \, \nabla \hat{u}_k}{(1+\eta'(\lambda_k(\hat{u}_{k-1}-u_k))) \, |\nabla \hat{u}_{k-1}| + (1-\eta'(\lambda_k(\hat{u}_{k-1}-u_k))) \, |\nabla u_k|} 
\end{align*} 
at the point $x$. The vector on the right hand side is a nonnegative multiple of $\hat{\nu}_k$. Using the identity 
\begin{align*} 
&|\sin(\alpha+\varphi) \, \hat{\nu}_{k-1} + \sin(\alpha-\varphi) \, \nu_k|^2 \\ 
&= \sin^2(\alpha+\varphi) + \sin^2(\alpha-\varphi) + 2 \sin(\alpha+\varphi) \sin(\alpha-\varphi) \cos(2\alpha) \\ 
&= \sin^2(2\alpha), 
\end{align*}
we conclude that 
\[\frac{\sin(\alpha+\varphi) \, \hat{\nu}_{k-1} + \sin(\alpha-\varphi) \, \nu_k}{\sin(2\alpha)} = \hat{\nu}_k\] 
at the point $x$. This completes the proof of Proposition \ref{recursive.relation.for.hat.nu}. \\

\begin{proposition}
\label{recursive.relation.for.hat.N}
If $\lambda_0$ is sufficiently large, we can find a collection of sets $W_0,\hdots,W_q \subset \mathbb{R}^n$ and a collection of maps $\hat{N}_0: W_0 \to S^{n-1},\hdots,\hat{N}_q: W_q \to S^{n-1}$ with the following properties: \\ 
(i) For each $0 \leq k \leq q$, $W_k$ is open and the map $\hat{N}_k: W_k \to S^{n-1}$ is smooth. \\
(ii) We have $W_0 = \mathbb{R}^n$. Moreover, the map $\hat{N}_0: W_0 \to S^{n-1}$ is constant and equal to $N_0$. \\
(iii) For each $1 \leq k \leq q$, we have $W_k = (W_{k-1} \setminus A_k) \cup \{\hat{u}_{k-1}-u_k < -\lambda_k^{-1}\}$, where 
\[A_k = W_{k-1} \cap \{|\hat{u}_{k-1} - u_k| \leq \lambda_k^{-1}\} \cap (\{d\hat{u}_{k-1} \wedge du_k = 0\} \cup \{\hat{N}_{k-1} \wedge N_k=0\}).\] 
(iv) For each $0 \leq k \leq q$, we have $\Omega \cap \{-2^{-k} \, \Xi^{-1} \leq \hat{u}_k \leq 0\} \subset W_k$. \\
(v) Let $1 \leq k \leq q$. Suppose that $x$ is a point in $W_k$ satisfying $\hat{u}_{k-1}(x)-u_k(x) > \lambda_k^{-1}$. Then $x \in W_{k-1}$ and $\hat{N}_k = \hat{N}_{k-1}$ at the point $x$. \\
(vi) Let $1 \leq k \leq q$. Suppose that $x$ is a point in $W_k$ satisfying $\hat{u}_{k-1}(x)-u_k(x) < -\lambda_k^{-1}$. Then $\hat{N}_k = N_k$ at the point $x$. \\
(vii) Let $1 \leq k \leq q$. Suppose that $x$ is a point in $W_k$ satisfying $|\hat{u}_{k-1}(x)-u_k(x)| \leq \lambda_k^{-1}$. Then $x \in W_{k-1} \setminus A_k$. We can find real numbers $\alpha \in (0,\frac{\pi}{2})$ and $\theta \in (0,\frac{\pi}{2\alpha})$ such that $\cos(2\alpha) = \langle \hat{\nu}_{k-1},\nu_k \rangle$ and $\cos(2\theta\alpha) = \langle \hat{N}_{k-1},N_k \rangle$ at the point $x$. Let us define a real number $\varphi \in [-\alpha,\alpha]$ by 
\begin{align*} 
&\frac{\tan(\varphi)}{\tan(\alpha)} \\ 
&= \frac{(1+\eta'(\lambda_k(\hat{u}_{k-1}-u_k))) \, |\nabla \hat{u}_{k-1}| - (1-\eta'(\lambda_k(\hat{u}_{k-1}-u_k))) \, |\nabla u_k|}{(1+\eta'(\lambda_k(\hat{u}_{k-1}-u_k))) \, |\nabla \hat{u}_{k-1}| + (1-\eta'(\lambda_k(\hat{u}_{k-1}-u_k))) \, |\nabla u_k|}, 
\end{align*} 
where $|\nabla \hat{u}_{k-1}|$ and $|\nabla u_k|$ are computed with respect to the metric $g$. Then 
\[\hat{N}_k = \frac{\sin(\theta(\alpha+\varphi)) \, \hat{N}_{k-1} + \sin(\theta(\alpha-\varphi)) \, N_k}{\sin(2\theta\alpha)}\] 
at the point $x$. \\
(viii) Let $0 \leq k \leq q$. Suppose that $x$ is a point in $\Omega \cap \{-2^{-k} \, \Xi^{-1} \leq \hat{u}_k \leq 0\}$. Then we can find nonnegative real numbers $a_0,\hdots,a_k$ such that $\hat{N}_k = \sum_{i=0}^k a_i \, N_i$ at the point $x$. Moreover, $a_i=0$ for all $0 \leq i \leq k$ satisfying $u_i(x) < \hat{u}_k(x) - 2k \lambda_i^{-1}$.
\end{proposition}

\textbf{Proof.} 
We proceed by induction. We define $W_0 = \mathbb{R}^n$. Moreover, we define the map $\hat{N}_0: W_0 \to S^{n-1}$ to be constant and equal to $N_0$. 

Suppose now that $1 \leq k \leq q$, and that we have constructed sets $W_0,\hdots,W_{k-1}$ and maps $\hat{N}_0: W_0 \to S^{n-1},\hdots,\hat{N}_{k-1}: W_{k-1} \to S^{n-1}$ satisfying properties (i)--(viii) above. We define 
\[A_k = W_{k-1} \cap \{|\hat{u}_{k-1} - u_k| \leq \lambda_k^{-1}\} \cap (\{d\hat{u}_{k-1} \wedge du_k = 0\} \cup \{\hat{N}_{k-1} \wedge N_k=0\})\] 
and 
\[W_k = (W_{k-1} \setminus A_k) \cup \{\hat{u}_{k-1}-u_k < -\lambda_k^{-1}\}.\] 
It follows from the induction hypothesis that $W_{k-1}$ is open. Since $A_k$ is a relatively closed subset of $W_{k-1}$, it follows that $W_k$ is open. 

It is clear that $W_k$ satisfies property (iii). Moreover, the induction hypothesis implies that $\Omega \cap \{-2^{-k+1} \, \Xi^{-1} \leq \hat{u}_{k-1} \leq 0\} \subset W_{k-1}$. Using Lemma \ref{inequality.for.hat.u_k}, we deduce that $\Omega \cap \{-2^{-k} \, \Xi^{-1} \leq \hat{u}_k \leq 0\} \subset W_{k-1} \cup \{\hat{u}_{k-1}-u_k < -\lambda_k^{-1}\}$. On the other hand, using Lemma \ref{inequality.for.hat.u_k} and Lemma \ref{transversality.2}, we obtain $d\hat{u}_{k-1} \wedge du_k \neq 0$ at each point in $\Omega \cap \{-2^{-k} \, \Xi^{-1} \leq \hat{u}_k \leq 0\} \cap \{|\hat{u}_{k-1} - u_k| \leq \lambda_k^{-1}\}$. Moreover, since the map $\hat{N}_{k-1}$ satisfies property (viii), it follows from Lemma \ref{transversality.1} and Lemma \ref{inequality.for.hat.u_k} that $\hat{N}_{k-1} \wedge N_k \neq 0$ at each point in $\Omega \cap \{-2^{-k} \, \Xi^{-1} \leq \hat{u}_k \leq 0\} \cap \{|\hat{u}_{k-1} - u_k| \leq \lambda_k^{-1}\}$. Therefore, the set $\Omega \cap \{-2^{-k} \, \Xi^{-1} \leq \hat{u}_k \leq 0\}$ is disjoint from $A_k$. Putting these facts together, we conclude that 
\[\Omega \cap \{-2^{-k} \, \Xi^{-1} \leq \hat{u}_k \leq 0\} \subset W_k.\] 
Therefore, $W_k$ satisfies property (iv).

In the next step, we define the map $\hat{N}_k: W_k \to S^{n-1}$. To that end, we distinguish three cases: 

\textit{Case 1:} Suppose that $x$ is a point in $W_k$ satisfying $\hat{u}_{k-1}(x)-u_k(x) > \lambda_k^{-1}$. This implies $x \in W_{k-1}$. In this case, we define $\hat{N}_k = \hat{N}_{k-1}$ at the point $x$.

\textit{Case 2:} Suppose that $x$ is a point in $W_k$ satisfying $\hat{u}_{k-1}(x)-u_k(x) < -\lambda_k^{-1}$. In this case, we define $\hat{N}_k = N_k$ at the point $x$.

\textit{Case 3:} Suppose that $x$ is a point in $W_k$ satisfying $|\hat{u}_{k-1}(x)-u_k(x)| \leq \lambda_k^{-1}$. Then $x \in W_{k-1} \setminus A_k$. This implies $d\hat{u}_{k-1} \wedge du_k \neq 0$, and $\hat{N}_{k-1} \wedge N_k \neq 0$ at the point $x$. Consequently, we can find real numbers $\alpha \in (0,\frac{\pi}{2})$ and $\theta \in (0,\frac{\pi}{2\alpha})$ such that $\cos(2\alpha) = \langle \hat{\nu}_{k-1},\nu_k \rangle$ and $\cos(2\theta\alpha) = \langle \hat{N}_{k-1},N_k \rangle$ at the point $x$. We define a real number $\varphi \in [-\alpha,\alpha]$ so that 
\begin{align*} 
&\frac{\tan(\varphi)}{\tan(\alpha)} \\ 
&= \frac{(1+\eta'(\lambda_k(\hat{u}_{k-1}-u_k))) \, |\nabla \hat{u}_{k-1}| - (1-\eta'(\lambda_k(\hat{u}_{k-1}-u_k))) \, |\nabla u_k|}{(1+\eta'(\lambda_k(\hat{u}_{k-1}-u_k))) \, |\nabla \hat{u}_{k-1}| + (1-\eta'(\lambda_k(\hat{u}_{k-1}-u_k))) \, |\nabla u_k|} 
\end{align*} 
at the point $x$. Since $d\hat{u}_{k-1} \wedge du_k \neq 0$ at the point $x$, the denominator is strictly positive; this ensures that $\varphi$ is well-defined. With this understood, we define 
\[\hat{N}_k = \frac{\sin(\theta(\alpha+\varphi)) \, \hat{N}_{k-1} + \sin(\theta(\alpha-\varphi)) \, N_k}{\sin(2\theta\alpha)}\] 
at the point $x$. In view of the identity 
\begin{align*} 
&|\sin(\theta(\alpha+\varphi)) \, \hat{N}_{k-1} + \sin(\theta(\alpha-\varphi)) \, N_k|^2 \\ 
&= \sin^2(\theta(\alpha+\varphi)) + \sin^2(\theta(\alpha-\varphi)) + 2 \sin(\theta(\alpha+\varphi)) \sin(\theta(\alpha-\varphi)) \cos(2\theta\alpha) \\ 
&= \sin^2(2\theta\alpha), 
\end{align*}
the map $\hat{N}_k$ takes values in $S^{n-1}$. 

It is clear from the definition that the map $\hat{N}_k: W_k \to S^{n-1}$ satisfies properties (v)--(vii). It remains to show that the map $\hat{N}_k$ satisfies property (viii). To prove this, we again distinguish three cases: 

\textit{Case 1:} Suppose that $x$ is a point in $\Omega$ satisfying $-2^{-k} \, \Xi^{-1} \leq \hat{u}_k(x) \leq 0$ and $\hat{u}_{k-1}(x)-u_k(x) > \lambda_k^{-1}$. In this case, $\hat{u}_k(x) = \hat{u}_{k-1}(x)$. In view of the induction hypothesis, we can find nonnegative real numbers $b_0,\hdots,b_{k-1}$ such that $\hat{N}_{k-1} = \sum_{i=0}^{k-1} b_i \, N_i$ at the point $x$. Moreover, $b_i=0$ for all $0 \leq i \leq k-1$ satisfying $u_i(x) < \hat{u}_{k-1}(x) - 2(k-1) \lambda_i^{-1}$. We define nonnegative real numbers $a_0,\hdots,a_k$ by $a_i = b_i$ for $0 \leq i \leq k-1$ and $a_k = 0$. Then $\hat{N}_k = \hat{N}_{k-1} = \sum_{i=0}^{k-1} b_i \, N_i = \sum_{i=0}^k a_i \, N_i$ at the point $x$. Moreover, $a_i=0$ for all $0 \leq i \leq k$ satisfying $u_i(x) < \hat{u}_k(x) - 2k \lambda_i^{-1}$. 

\textit{Case 2:} Suppose that $x$ is a point in $\Omega$ satisfying $-2^{-k} \, \Xi^{-1} \leq \hat{u}_k(x) \leq 0$ and $\hat{u}_{k-1}(x)-u_k(x) < -\lambda_k^{-1}$. In this case, $\hat{u}_k(x) = u_k(x)$. We define nonnegative real numbers $a_0,\hdots,a_k$ by $a_i=0$ for $0 \leq i \leq k-1$ and $a_k=1$. Then $\hat{N}_k = N_k = \sum_{i=0}^k a_i \, N_i$. Moreover, $a_i=0$ for all $0 \leq i \leq k$ satisfying $u_i(x) < \hat{u}_k(x) - 2k \lambda_i^{-1}$. 

\textit{Case 3:} Suppose that $x$ is a point in $\Omega$ satisfying $-2^{-k} \, \Xi^{-1} \leq \hat{u}_k(x) \leq 0$ and $|\hat{u}_{k-1}(x)-u_k(x)| \leq \lambda_k^{-1}$. Using Lemma \ref{inequality.for.hat.u_k}, we obtain 
\[\hat{u}_k(x) \leq \max \{\hat{u}_{k-1}(x),u_k(x)\} + \lambda_k^{-1} \leq \min \{\hat{u}_{k-1}(x),u_k(x)\} + 2\lambda_k^{-1}.\] 
Let $\alpha \in (0,\frac{\pi}{2})$, $\theta \in (0,\frac{\pi}{2\alpha})$, and $\varphi \in [-\alpha,\alpha]$ be defined as above. In view of the induction hypothesis, we can find nonnegative real numbers $b_0,\hdots,b_{k-1}$ such that $\hat{N}_{k-1} = \sum_{i=0}^{k-1} b_i \, N_i$ at the point $x$. Moreover, $b_i=0$ for all $0 \leq i \leq k-1$ satisfying $u_i(x) < \hat{u}_{k-1}(x) - 2(k-1) \lambda_i^{-1}$. We define 
\[a_i = \frac{\sin(\theta(\alpha+\varphi))}{\sin(2\theta\alpha)} \, b_i\] 
for $0 \leq i \leq k-1$ and 
\[a_k = \frac{\sin(\theta(\alpha-\varphi))}{\sin(2\theta\alpha)}.\] 
Since $\theta\alpha \in (0,\frac{\pi}{2})$ and $\varphi \in [-\alpha,\alpha]$, it follows that $a_0,\hdots,a_k$ are nonnegative real numbers. Moreover, 
\begin{align*} 
\hat{N}_k 
&= \frac{\sin(\theta(\alpha+\varphi)) \, \hat{N}_{k-1} + \sin(\theta(\alpha-\varphi)) \, N_k}{\sin(2\theta\alpha)} \\ 
&= \frac{\sin(\theta(\alpha+\varphi)) \, \sum_{i=0}^{k-1} b_i \, N_i + \sin(\theta(\alpha-\varphi)) \, N_k}{\sin(2\theta\alpha)} \\ 
&= \sum_{i=0}^k a_i \, N_i 
\end{align*}
at the point $x$. Finally, $a_i=0$ for all $0 \leq i \leq k$ satisfying $u_i(x) < \hat{u}_k(x) - 2k \lambda_i^{-1}$. 

To summarize, we have shown that the map $\hat{N}_k$ satisfies properties (v)--(viii). This completes the proof of Proposition \ref{recursive.relation.for.hat.N}. \\

\begin{remark} 
Let us explain the geometric significance of Proposition \ref{recursive.relation.for.hat.nu} (iii) and Proposition \ref{recursive.relation.for.hat.N} (vii). Let $1 \leq k \leq q$. Suppose that $x$ is a point in $\Omega$ with the property that $2^{-k} \, \Xi^{-1} \leq \hat{u}_k(x) \leq 0$ and $|\hat{u}_{k-1}(x)-u_k(x)| \leq \lambda_k^{-1}$. Let $\alpha \in (0,\frac{\pi}{2})$, $\theta \in (0,\frac{\pi}{2\alpha})$, and $\varphi \in [-\alpha,\alpha]$ be defined as above. Then $\varangle (\hat{\nu}_{k-1},\nu_k) = 2\alpha$, $\varangle (\hat{\nu}_{k-1},\hat{\nu}_k) = \alpha-\varphi$, and $\varangle (\nu_k,\hat{\nu}_k) = \alpha+\varphi$. Moreover, $\varangle (\hat{N}_{k-1},N_k) = 2\theta\alpha$, $\varangle (\hat{N}_{k-1},\hat{N}_k) = \theta(\alpha-\varphi)$, and $\varangle (N_k,\hat{N}_k) = \theta(\alpha+\varphi)$.
\end{remark}

\begin{proposition}
\label{transversality.3}
If $\lambda_0$ is sufficiently large, then the following holds. Let $0 \leq j < k \leq q$. If $x$ is a point in $\Omega$ satisfying $-2^{-j} \, \Xi^{-1} \leq \hat{u}_j(x) \leq 0$ and $-\Xi^{-1} \leq u_k(x) \leq 0$, then $|\hat{N}_j \wedge N_k| \geq \Xi^{-1}$.
\end{proposition}

\textbf{Proof.} 
Consider a point $x \in \Omega$ satisfying $-2^{-j} \, \Xi^{-1} \leq \hat{u}_j(x) \leq 0$ and $-\Xi^{-1} \leq u_k(x) \leq 0$. By Proposition \ref{recursive.relation.for.hat.N} (viii), we can find nonnegative real numbers $a_0,\hdots,a_j$ such that $\hat{N}_j = \sum_{i=0}^j a_i \, N_i$ at the point $x$. Moreover, $a_i=0$ for all $i \in \{0,1,\hdots,j\}$ satisfying $u_i(x) < \hat{u}_j(x) - 2j \lambda_i^{-1}$. Since $-2^{-j} \, \Xi^{-1} \leq \hat{u}_j(x) \leq 0$, it follows that $a_i=0$ for all $i \in \{0,1,\hdots,j\}$ satisfying $u_i(x) \leq -2\Xi^{-1}$. Using Lemma \ref{transversality.1}, we obtain 
\[1 = |\hat{N}_j| \leq \sum_{i=0}^j a_i \leq \Xi \, \Big | \Big ( \sum_{i=0}^j a_i \, N_i \Big ) \wedge N_k \Big | = \Xi \, |\hat{N}_j \wedge N_k|\] 
at the point $x$. This completes the proof of Proposition \ref{transversality.3}. \\

\begin{proposition}
\label{C1.bound.for.hat.N}
Let $0 \leq k \leq q$. Then $|d\hat{N}_k| \leq C\lambda_k$ at each point in $\Omega \cap \{-2^{-k} \, \Xi^{-1} \leq \hat{u}_k \leq 0\}$, where $C$ is independent of $\gamma$ and $\lambda_0$.
\end{proposition}

\textbf{Proof.} 
We argue by induction on $k$. Since the map $\hat{N}_0$ is constant, the assertion is obvious for $k=0$.

Suppose now that $1 \leq k \leq q$, and that $|d\hat{N}_{k-1}| \leq C\lambda_{k-1}$ at each point in $\Omega \cap \{-2^{-k+1} \, \Xi^{-1} \leq \hat{u}_{k-1} \leq 0\}$, where $C$ is independent of $\gamma$ and $\lambda_0$. Let us consider an arbitrary point $x_0 \in \Omega \cap \{-2^{-k} \, \Xi^{-1} \leq \hat{u}_k \leq 0\}$. We distinguish three cases:

\textit{Case 1:} Suppose that $\hat{u}_{k-1}(x_0)-u_k(x_0) > \lambda_k^{-1}$. In this case, $\hat{N}_k = \hat{N}_{k-1}$ in an open neighborhood of the point $x_0$. Therefore, $|d\hat{N}_k| = |d\hat{N}_{k-1}| \leq C\lambda_k$ at $x_0$.

\textit{Case 2:} Suppose that $\hat{u}_{k-1}(x_0)-u_k(x_0) < -\lambda_k^{-1}$. In this case, $\hat{N}_k = N_k$ in an open neighborhood of the point $x_0$. Therefore, $d\hat{N}_k = 0$ at $x_0$.

\textit{Case 3:} Suppose that $|\hat{u}_{k-1}(x_0)-u_k(x_0)| \leq \lambda_k^{-1}$. In this case, $\alpha$, $\theta$, and $\varphi$ can be viewed as smooth functions which are defined in an open neighborhood of the point $x_0$. Differentiating the identity $\cos(2\alpha) = \langle \hat{\nu}_{k-1},\nu_k \rangle$, we obtain $|d\cos(2\alpha)| \leq C\lambda_{k-1}$ at $x_0$. By Lemma \ref{transversality.2}, $\alpha \in (0,\frac{\pi}{2})$ is bounded away from $0$ and $\frac{\pi}{2}$ at the point $x_0$. Consequently, $|d\alpha| \leq C\lambda_{k-1}$ at $x_0$. Differentiating the identity $\cos(2\theta\alpha) = \langle \hat{N}_{k-1},N_k \rangle$, we obtain $|d\cos(2\theta\alpha)| \leq C\lambda_{k-1}$ at $x_0$. By Proposition \ref{transversality.3}, $\theta\alpha \in (0,\frac{\pi}{2})$ is bounded away from $0$ and $\frac{\pi}{2}$ at the point $x_0$. Consequently, $|d(\theta\alpha)| \leq C\lambda_{k-1}$ at $x_0$. This implies $|d\theta| \leq C\lambda_{k-1}$ at $x_0$. A similar argument gives $|d\varphi| \leq C\lambda_k$ at $x_0$. Since 
\[\hat{N}_k = \frac{\sin(\theta(\alpha+\varphi)) \, \hat{N}_{k-1} + \sin(\theta(\alpha-\varphi)) \, N_k}{\sin(2\theta\alpha)}\] 
in an open neighborhood of $x_0$, we conclude that $|d\hat{N}_k| \leq C\lambda_k$ at $x_0$. This completes the proof of Proposition \ref{C1.bound.for.hat.N}. \\

\begin{proposition}
\label{hat.N.Euclidean.case} 
Suppose that $g$ is the Euclidean metric. For each $0 \leq k \leq q$ and each point $x \in \Omega \cap \{-2^{-k} \, \Xi^{-1} \leq \hat{u}_k \leq 0\}$, the vector $\hat{N}_k(x)$ is the unit normal vector to the level set of $\hat{u}_k$ passing through $x$.
\end{proposition} 

\textbf{Proof.} 
This can be proved by induction on $k$, using the recursive relations in Proposition \ref{recursive.relation.for.hat.nu} and Proposition \ref{recursive.relation.for.hat.N}. \\

\begin{definition} 
We denote by $\hat{N}$ the restriction of the map $\hat{N}_q$ to the hypersurface $\hat{\Sigma} = \{\hat{u}_q=0\}$. 
\end{definition}

\begin{proposition} 
The map $\hat{N}: \hat{\Sigma} \to S^{n-1}$ is homotopic to the Gauss map of $\hat{\Sigma}$.
\end{proposition}

\textbf{Proof.} 
In the special case when $g$ is the Euclidean metric, the assertion follows from Proposition \ref{hat.N.Euclidean.case}. To prove the assertion in general, we deform the metric $g$ to the Euclidean metric. \\

\section{Controlling the angles under the smoothing procedure}

\label{controlling.angles}

Suppose that $n \geq 3$ is an integer, and $\Omega = \bigcap_{m=0}^q \{u_m \leq 0\}$ is a compact, convex polytope in $\mathbb{R}^n$ with non-empty interior. Throughout this section, we assume that Assumptions \ref{no.redundant.inequalities}, \ref{gradient.of.u_k}, and \ref{angles.bounded.by.pi/2} are satisfied. Moreover, we assume that $g$ is a Riemannian metric on $\mathbb{R}^n$ which satisfies the following assumption. 

\begin{assumption}
\label{angle.comparison}
Let $0 \leq j < k \leq q$. If $x$ is a point in $\Omega$ with $u_j(x)=u_k(x)=0$, then $\langle \nu_j,\nu_k \rangle \leq \langle N_j,N_k \rangle$ at the point $x$.
\end{assumption}

\begin{lemma}
\label{approximate.angle.inequality.2}
Given $\varepsilon \in (0,1)$, we can find a small positive real number $\delta$ (depending on $\varepsilon$) with the following property. Let $0 \leq j < k \leq q$. If $x$ is a point in $\Omega$ with $-\delta \leq u_j(x) \leq 0$ and $-\delta \leq u_k(x) \leq 0$, then $\langle \nu_j,\nu_k \rangle - \varepsilon \leq \langle N_j,N_k \rangle$ at the point $x$.
\end{lemma}

\textbf{Proof.} 
Suppose that the assertion is false. We consider a sequence of counterexamples and pass to the limit. In the limit, we obtain a pair of integers $0 \leq j < k \leq q$ and a point $x \in \Omega$ such that $u_j(x)=u_k(x)= 0$ and $\langle \nu_j,\nu_k \rangle - \varepsilon \geq \langle N_j,N_k \rangle$ at the point $x$. This contradicts Assumption \ref{angle.comparison}. This completes the proof of Lemma \ref{approximate.angle.inequality.2}. \\

\begin{lemma}
\label{preservation.of.angle.inequality.under.smoothing}
We can find a constant $Q>1$ (independent of $\gamma$ and $\lambda_0$) such that the following holds. Consider a pair of integers $1 \leq j < k \leq q$. Moreover, suppose that $x$ is a point in $\Omega$ satisfying $-2^{-j} \, \Xi^{-1} \leq \hat{u}_j(x) \leq 0$ and $|\hat{u}_{j-1}(x) - u_j(x)| \leq \lambda_j^{-1}$. If 
\begin{align*} 
&\langle \hat{\nu}_{j-1},\nu_j \rangle - \varepsilon \leq \langle \hat{N}_{j-1},N_j \rangle, \\ 
&\langle \hat{\nu}_{j-1},\nu_k \rangle - \varepsilon \leq \langle \hat{N}_{j-1},N_k \rangle \leq \varepsilon, \\ 
&\langle \nu_j,\nu_k \rangle - \varepsilon \leq \langle N_j,N_k \rangle \leq \varepsilon 
\end{align*} 
at the point $x$, then 
\[\langle \hat{\nu}_j,\nu_k \rangle - Q\varepsilon \leq \langle \hat{N}_j,N_k \rangle \leq Q\varepsilon\]
at the point $x$.
\end{lemma} 

\textbf{Proof.} 
We can find real numbers $\alpha \in (0,\frac{\pi}{2})$ and $\theta \in (0,\frac{\pi}{2\alpha})$ such that $\cos(2\alpha) = \langle \hat{\nu}_{j-1},\nu_j \rangle$ and $\cos(2\theta\alpha) = \langle \hat{N}_{j-1},N_j \rangle$ at the point $x$. As above, we define a real number $\varphi \in [-\alpha,\alpha]$ by 
\begin{align*} 
&\frac{\tan(\varphi)}{\tan(\alpha)} \\ 
&= \frac{(1+\eta'(\lambda_j(\hat{u}_{j-1}-u_j))) \, |\nabla \hat{u}_{j-1}| - (1-\eta'(\lambda_j(\hat{u}_{j-1}-u_j))) \, |\nabla u_j|}{(1+\eta'(\lambda_j(\hat{u}_{j-1}-u_j))) \, |\nabla \hat{u}_{j-1}| + (1-\eta'(\lambda_j(\hat{u}_{j-1}-u_j))) \, |\nabla u_j|}.
\end{align*} 
Using Proposition \ref{recursive.relation.for.hat.nu} (iii), we obtain 
\[\hat{\nu}_j = \frac{\sin(\alpha+\varphi) \, \hat{\nu}_{j-1} + \sin(\alpha-\varphi) \, \nu_j}{\sin(2\alpha)}\] 
at the point $x$. Moreover, Proposition \ref{recursive.relation.for.hat.N} (vii) implies that 
\[\hat{N}_j = \frac{\sin(\theta(\alpha+\varphi)) \, \hat{N}_{j-1} + \sin(\theta(\alpha-\varphi)) \, N_j}{\sin(2\theta\alpha)}\] 
at the point $x$. This gives
\begin{align} 
\label{id.1}
\frac{\cos(\theta\varphi)}{\cos(\theta\alpha)} \, \varepsilon - \langle \hat{N}_j,N_k \rangle \notag
&= \frac{\sin(\theta(\alpha+\varphi))}{\sin(2\theta\alpha)} \, (\varepsilon - \langle \hat{N}_{j-1},N_k \rangle) \\ 
&+ \frac{\sin(\theta(\alpha-\varphi))}{\sin(2\theta\alpha)} \, (\varepsilon - \langle N_j,N_k \rangle) 
\end{align} 
and 
\begin{align} 
\label{id.2}
&\langle \hat{N}_j,N_k \rangle - \langle \hat{\nu}_j,\nu_k \rangle + 2 \, \frac{\cos(\varphi)}{\cos(\alpha)} \, \varepsilon - \frac{\cos(\theta \varphi)}{\cos(\theta\alpha)} \, \varepsilon \notag \\ 
&= \Big ( \frac{\sin(\alpha+\varphi)}{\sin(2\alpha)} - \frac{\sin(\theta(\alpha+\varphi))}{\sin(2\theta\alpha)} \Big ) \, (\varepsilon - \langle \hat{N}_{j-1},N_k \rangle) \notag \\ 
&+ \Big ( \frac{\sin(\alpha-\varphi)}{\sin(2\alpha)} - \frac{\sin(\theta(\alpha-\varphi))}{\sin(2\theta\alpha)} \Big ) \, (\varepsilon - \langle N_j,N_k \rangle) \\ 
&+ \frac{\sin(\alpha+\varphi)}{\sin(2\alpha)} \, (\langle \hat{N}_{j-1},N_k \rangle - \langle \hat{\nu}_{j-1},\nu_k \rangle + \varepsilon) \notag \\ 
&+ \frac{\sin(\alpha-\varphi)}{\sin(2\alpha)} \, (\langle N_j,N_k \rangle - \langle \nu_j,\nu_k \rangle + \varepsilon) \notag
\end{align} 
at the point $x$. By assumption, 
\begin{align*} 
&\varepsilon - \langle \hat{N}_{j-1},N_k \rangle \geq 0, \\ 
&\varepsilon - \langle N_j,N_k \rangle \geq 0, \\
&\langle \hat{N}_{j-1},N_k \rangle - \langle \hat{\nu}_{j-1},\nu_k \rangle + \varepsilon \geq 0, \\ 
&\langle N_j,N_k \rangle - \langle \nu_j,\nu_k \rangle + \varepsilon \geq 0 
\end{align*} 
at the point $x$. Using the identity (\ref{id.1}), we obtain 
\begin{equation} 
\label{ineq.1}
\frac{\cos(\theta\varphi)}{\cos(\theta\alpha)} \, \varepsilon - \langle \hat{N}_j,N_k \rangle \geq 0 
\end{equation}
at the point $x$. 

In the next step, we bound $\langle \hat{N}_j,N_k \rangle - \langle \hat{\nu}_j,\nu_k \rangle$ from below. To that end, we distinguish two cases:

\textit{Case 1:} Suppose that $\theta \in (0,1)$. Applying Proposition \ref{elementary.inequality.2} with $\beta = \frac{\alpha+\varphi}{2} \in [0,\alpha]$ gives 
\[\frac{\sin(\alpha+\varphi)}{\sin(2\alpha)} - \frac{\sin(\theta(\alpha+\varphi))}{\sin(2\theta\alpha)} \geq 0.\] 
Similarly, applying Proposition \ref{elementary.inequality.2} with $\beta = \frac{\alpha-\varphi}{2} \in [0,\alpha]$ gives 
\[\frac{\sin(\alpha-\varphi)}{\sin(2\alpha)} - \frac{\sin(\theta(\alpha-\varphi))}{\sin(2\theta\alpha)} \geq 0.\]
Using the identity (\ref{id.2}), we conclude that 
\begin{equation} 
\label{ineq.2}
\langle \hat{N}_j,N_k \rangle - \langle \hat{\nu}_j,\nu_k \rangle + 2 \, \frac{\cos(\varphi)}{\cos(\alpha)} \, \varepsilon - \frac{\cos(\theta\varphi)}{\cos(\theta\alpha)} \, \varepsilon \geq 0 
\end{equation}
at the point $x$. 

\textit{Case 2:} Suppose that $\theta \in [1,\frac{\pi}{2\alpha})$. Applying Proposition \ref{elementary.inequality.3} with $\beta = \frac{\alpha+\varphi}{2} \in [0,\alpha]$ gives 
\[\frac{\sin(\alpha+\varphi)}{\sin(2\alpha)} - \frac{\sin(\theta(\alpha+\varphi))}{\sin(2\theta\alpha)} \geq -\frac{4(\theta-1)\alpha}{\sin(2\alpha) \sin(2\theta\alpha)}.\] 
Similarly, applying Proposition \ref{elementary.inequality.3} with $\beta = \frac{\alpha-\varphi}{2} \in [0,\alpha]$ gives 
\[\frac{\sin(\alpha-\varphi)}{\sin(2\alpha)} - \frac{\sin(\theta(\alpha-\varphi))}{\sin(2\theta\alpha)} \geq -\frac{4(\theta-1)\alpha}{\sin(2\alpha) \sin(2\theta\alpha)}.\] 
Using the identity (\ref{id.2}), we obtain 
\begin{align*} 
&\langle \hat{N}_j,N_k \rangle - \langle \hat{\nu}_j,\nu_k \rangle + 2 \, \frac{\cos(\varphi)}{\cos(\alpha)} \, \varepsilon - \frac{\cos(\theta\varphi)}{\cos(\theta\alpha)} \, \varepsilon \\ 
&\geq -\frac{4(\theta-1)\alpha}{\sin(2\alpha) \sin(2\theta\alpha)} \, (2\varepsilon - \langle \hat{N}_{j-1},N_k \rangle - \langle N_j,N_k \rangle) \\ 
&\geq -\frac{8(\theta-1)\alpha (1+\varepsilon)}{\sin(2\alpha) \sin(2\theta\alpha)} 
\end{align*}
at the point $x$. On the other hand, 
\[2(\theta-1)\alpha \inf_{t \in [2\alpha,2\theta\alpha]} \sin(t) \leq \cos(2\alpha) - \cos(2\theta\alpha) = \langle \hat{\nu}_{j-1},\nu_j \rangle - \langle \hat{N}_{j-1},N_j \rangle \leq \varepsilon\] 
at the point $x$. Putting these facts together, we conclude that 
\begin{align} 
\label{ineq.3}
&\langle \hat{N}_j,N_k \rangle - \langle \hat{\nu}_j,\nu_k \rangle + 2 \, \frac{\cos(\varphi)}{\cos(\alpha)} \, \varepsilon - \frac{\cos(\theta\varphi)}{\cos(\theta\alpha)} \, \varepsilon \notag \\ 
&\geq -\frac{4\varepsilon(1+\varepsilon)}{\sin(2\alpha) \sin(2\theta\alpha)  \inf_{t \in [2\alpha,2\theta\alpha]} \sin(t)} 
\end{align}
at the point $x$. 

By Lemma \ref{transversality.2}, $\alpha \in (0,\frac{\pi}{2})$ is bounded away from $0$ and $\frac{\pi}{2}$ at the point $x$. By Proposition \ref{transversality.3}, $\theta\alpha \in (0,\frac{\pi}{2})$ is bounded away from $0$ and $\frac{\pi}{2}$ at the point $x$. Combining (\ref{ineq.1}), (\ref{ineq.2}), and (\ref{ineq.3}), the assertion follows. This completes the proof of Lemma \ref{preservation.of.angle.inequality.under.smoothing}. \\

\begin{proposition} 
\label{approximate.angle.inequality.3}
Let $0 \leq j \leq q-1$. Given $\varepsilon \in (0,1)$, we can find a real number $\delta \in (0,2^{-j} \, \Xi^{-1})$ (depending on $\varepsilon$, but independent of $\gamma$ and $\lambda_0$) with the following property. Suppose that $\lambda_0 \geq 2\delta^{-1}$, $k \in \{j+1,\hdots,q\}$, and $x$ is a point in $\Omega$ satisfying $-\delta \leq \hat{u}_j(x) \leq 0$ and $-\delta \leq u_k(x) \leq 0$. Then $\langle \hat{\nu}_j,\nu_k \rangle - \varepsilon \leq \langle \hat{N}_j,N_k \rangle \leq \varepsilon$ at the point $x$.
\end{proposition}

\textbf{Proof.} We argue by induction on $j$. For $j=0$, the assertion follows from Lemma \ref{approximate.angle.inequality.1} and Lemma \ref{approximate.angle.inequality.2}. We next assume that $1 \leq j \leq q$, and that the assertion is true for $j-1$. We claim that the assertion is true for $j$. To prove this, let us fix a constant $Q>1$ so that the conclusion of Lemma \ref{preservation.of.angle.inequality.under.smoothing} holds. Note that $Q$ is independent of $\gamma$ and $\lambda_0$. Let $\varepsilon \in (0,1)$ be given. In view of the induction hypothesis, we can find a real number $\delta_0 \in (0,2^{-j+1} \, \Xi^{-1})$ (depending on $\varepsilon$, but independent of $\gamma$ and $\lambda_0$) such that 
\begin{equation} 
\label{property.of.delta_0}
\langle \hat{\nu}_{j-1},\nu_k \rangle - Q^{-1} \varepsilon \leq \langle \hat{N}_{j-1},N_k \rangle \leq Q^{-1} \varepsilon 
\end{equation}
whenever $\lambda_0 \geq 2\delta_0^{-1}$, $k \in \{j,\hdots,q\}$, and $x$ is a point in $\Omega$ satisfying $-\delta_0 \leq \hat{u}_{j-1}(x) \leq 0$ and $-\delta_0 \leq u_k(x) \leq 0$. By Lemma \ref{approximate.angle.inequality.1} and Lemma \ref{approximate.angle.inequality.2}, we can find a positive real number $\delta_1$ (depending on $\varepsilon$) such that 
\begin{equation} 
\label{property.of.delta_1}
\langle \nu_j,\nu_k \rangle - Q^{-1} \varepsilon \leq \langle N_j,N_k \rangle \leq Q^{-1} \varepsilon 
\end{equation}
whenever $k \in \{j+1,\hdots,q\}$ and $x$ is a point in $\Omega$ satisfying $-\delta_1 \leq u_j(x) \leq 0$ and $-\delta_1 \leq u_k(x) \leq 0$. 

We define $\delta := \frac{1}{2} \min \{\delta_0,\delta_1\} \in (0,2^{-j} \, \Xi^{-1})$. Note that $\delta$ may depend on $\varepsilon$, but is independent of $\gamma$ and $\lambda_0$. We claim that $\delta$ has the desired property. To prove this, we assume that $\lambda_0 \geq 2\delta^{-1}$. Moreover, we consider an integer $k \in \{j+1,\hdots,q\}$ and a point $x \in \Omega$ such that $-\delta \leq \hat{u}_j(x) \leq 0$ and $-\delta \leq u_k(x) \leq 0$. We distinguish three cases: 

\textit{Case 1:} Suppose first that $\hat{u}_{j-1}(x) - u_j(x) > \lambda_j^{-1}$. In this case, $\hat{u}_j(x) = \hat{u}_{j-1}(x)$. Consequently, $-\delta \leq \hat{u}_{j-1}(x) \leq 0$ and $-\delta \leq u_k(x) \leq 0$. Using (\ref{property.of.delta_0}), we obtain 
\[\langle \hat{\nu}_{j-1},\nu_k \rangle - Q^{-1} \varepsilon \leq \langle \hat{N}_{j-1},N_k \rangle \leq Q^{-1} \varepsilon\] 
at the point $x$. Since $\hat{\nu}_j = \hat{\nu}_{j-1}$ and $\hat{N}_j = \hat{N}_{j-1}$ at the point $x$, it follows that 
\[\langle \hat{\nu}_j,\nu_k \rangle - Q^{-1} \varepsilon \leq \langle \hat{N}_j,N_k \rangle \leq Q^{-1} \varepsilon\] 
at the point $x$.

\textit{Case 2:} Suppose next that $\hat{u}_{j-1}(x) - u_j(x) < -\lambda_j^{-1}$. In this case, $\hat{u}_j(x) = u_j(x)$. Consequently, $-\delta \leq u_j(x) \leq 0$ and $-\delta \leq u_k(x) \leq 0$. Using (\ref{property.of.delta_1}), we obtain 
\[\langle \nu_j,\nu_k \rangle - Q^{-1} \varepsilon \leq \langle N_j,N_k \rangle \leq Q^{-1} \varepsilon.\] 
Since $\hat{\nu}_j = \nu_j$ and $\hat{N}_j = N_j$ at the point $x$, it follows that 
\[\langle \hat{\nu}_j,\nu_k \rangle - Q^{-1} \varepsilon \leq \langle \hat{N}_j,N_k \rangle \leq Q^{-1} \varepsilon\] 
at the point $x$.

\textit{Case 3:} Suppose finally that $|\hat{u}_{j-1}(x) - u_j(x)| \leq \lambda_j^{-1}$. Since $\hat{u}_j(x) \in [-\delta,0]$, Lemma \ref{inequality.for.hat.u_k} implies that $\hat{u}_{j-1}(x) \in [-\delta-2\lambda_j^{-1},0]$ and $u_j(x) \in [-\delta-2\lambda_j^{-1},0]$. Since $\lambda_j \geq \lambda_0 \geq 2\delta^{-1}$, we conclude that $-2\delta \leq \hat{u}_{j-1}(x) \leq 0$ and $-2\delta \leq u_j(x) \leq 0$. Moreover, $-\delta \leq u_k(x) \leq 0$. Using (\ref{property.of.delta_0}), we obtain 
\[\langle \hat{\nu}_{j-1},\nu_j \rangle - Q^{-1} \varepsilon \leq \langle \hat{N}_{j-1},N_j \rangle \leq Q^{-1} \varepsilon\] 
and 
\[\langle \hat{\nu}_{j-1},\nu_k \rangle - Q^{-1} \varepsilon \leq \langle \hat{N}_{j-1},N_k \rangle \leq Q^{-1} \varepsilon\] 
at the point $x$. Using (\ref{property.of.delta_1}), we obtain 
\[\langle \nu_j,\nu_k \rangle - Q^{-1} \varepsilon \leq \langle N_j,N_k \rangle \leq Q^{-1} \varepsilon\] 
at the point $x$. Using Lemma \ref{preservation.of.angle.inequality.under.smoothing}, we conclude that 
\[\langle \hat{\nu}_j,\nu_k \rangle - \varepsilon \leq \langle \hat{N}_j,N_k \rangle \leq \varepsilon\] 
at the point $x$. This completes the proof of Proposition \ref{approximate.angle.inequality.3}. \\

\section{Estimating $\max \{\|d\hat{N}\|_{\text{\rm tr}} - H,0\}$}

\label{controlling.mean.curvature}

Suppose that $n \geq 3$ is an integer, and $\Omega = \bigcap_{m=0}^q \{u_m \leq 0\}$ is a compact, convex polytope in $\mathbb{R}^n$ with non-empty interior. Throughout this section, we assume that Assumptions \ref{no.redundant.inequalities}, \ref{gradient.of.u_k}, and \ref{angles.bounded.by.pi/2} are satisfied. Moreover, we assume that $g$ is a Riemannian metric on $\mathbb{R}^n$ which satisfies Assumption \ref{angle.comparison}. We further assume that the metric $g$ satisfies the following assumption.

\begin{assumption}
\label{faces.are.mean.convex}
For each $k \in \{0,1,\hdots,q\}$, the mean curvature of the hypersurface $\{u_k=0\}$ with respect to $g$ is nonnegative at each point in $\Omega \cap \{u_k=0\}$.
\end{assumption} 

In the remainder of this section, we will estimate the quantity $\max \{\|d\hat{N}\|_{\text{\rm tr}} - H,0\}$, where $H$ denotes the mean curvature of $\hat{\Sigma}$ with respect to the metric $g$.

\begin{proposition} 
\label{pointwise.bound.on.F_k}
Let $0 \leq k \leq q$. Then $H - \|d\hat{N}\|_{\text{\rm tr}} \geq 0$ on the set $F_k$.
\end{proposition}

\textbf{Proof.} Let us fix a point $x_0 \in F_k$. We consider a small open neighborhood $W$ of $x_0$. We choose $W$ so that 
\[W \subset \bigcap_{m=k+1}^q \{\hat{u}_{m-1}-u_m > \lambda_m^{-1}\}.\] 
If $k \neq 0$, we further require that $W \subset \{\hat{u}_{k-1}-u_k < -\lambda_k^{-1}\}$. 

Then $u_k = \hat{u}_k = \hat{u}_q$ at each point in $W$. This implies $\nu_k = \hat{\nu}_k = \hat{\nu}_q = \hat{\nu}$ at each point in $\hat{\Sigma} \cap W$. Moreover, it follows from statements (v) and (vi) in Proposition \ref{recursive.relation.for.hat.N} that $N_k = \hat{N}_k = \hat{N}_q = \hat{N}$ at each point in $\hat{\Sigma} \cap W$. Using Assumption \ref{faces.are.mean.convex}, we conclude that $H \geq 0$ and $\|d\hat{N}\|_{\text{\rm tr}} = 0$ at $x_0$. This completes the proof of Proposition \ref{pointwise.bound.on.F_k}. \\

\begin{proposition}
\label{pointwise.bound.on.E_jk}
Let $0 \leq j < k \leq q$. If $\lambda_0$ is sufficiently large, then $H - \|d\hat{N}\|_{\text{\rm tr}} \geq -C\lambda_k \max \{\langle \nu_j,\nu_k \rangle - \langle N_j,N_k \rangle,0\} - C$ on the set $E_{j,k}$. Here, $C$ is independent of $\gamma$ and $\lambda_0$. 
\end{proposition}

\textbf{Proof.} 
Let us fix a point $x_0 \in E_{j,k}$. In the following, we will consider a small open neighborhood $W$ of $x_0$. We choose $W$ so that 
\[W \subset \bigcap_{m=k+1}^q \{\hat{u}_{m-1}-u_m > \lambda_m^{-1}\} \cap \bigcap_{m=j+1}^{k-1} \{\hat{u}_{m-1}-u_m > \lambda_m^{-1}\}.\] 
If $j \neq 0$, we further require that $W \subset \{\hat{u}_{j-1}-u_j < -\lambda_j^{-1}\}$. 

Then $u_j = \hat{u}_j = \hat{u}_{k-1}$ and $\hat{u}_k = \hat{u}_q$ at each point in $W$. This implies $\nu_j = \hat{\nu}_j = \hat{\nu}_{k-1}$ and $\hat{\nu}_k = \hat{\nu}_q = \hat{\nu}$ at each point in $\hat{\Sigma} \cap W$. Moreover, it follows from statements (v) and (vi) in Proposition \ref{recursive.relation.for.hat.N} that $N_j = \hat{N}_j = \hat{N}_{k-1}$ and $\hat{N}_k = \hat{N}_q = \hat{N}$ at each point in $\hat{\Sigma} \cap W$. 

We define a smooth function $\alpha: \hat{\Sigma} \cap W \to (0,\frac{\pi}{2})$ by 
\begin{equation} 
\label{formula.for.alpha}
\cos(2\alpha) = \langle \nu_j,\nu_k \rangle 
\end{equation} 
at each point in $\hat{\Sigma} \cap W$. Moreover, we define a smooth function $\theta: \hat{\Sigma} \cap W \to \mathbb{R}$ so that $\theta \in (0,\frac{\pi}{2\alpha})$ at each point in $\hat{\Sigma} \cap W$ and 
\begin{equation} 
\label{formula.for.theta}
\cos(2\theta\alpha) = \langle N_j,N_k \rangle 
\end{equation} 
at each point in $\hat{\Sigma} \cap W$. Finally, we define a function $\varphi: \hat{\Sigma} \cap W \to \mathbb{R}$ so that $\varphi \in [-\alpha,\alpha]$ at each point in $\hat{\Sigma} \cap W$ and 
\begin{align} 
\label{formula.for.varphi}
&\frac{\tan(\varphi)}{\tan(\alpha)} \notag \\ 
&= \frac{(1+\eta'(\lambda_k(u_j-u_k))) \, |\nabla u_j| - (1-\eta'(\lambda_k(u_j-u_k))) \, |\nabla u_k|}{(1+\eta'(\lambda_k(u_j-u_k))) \, |\nabla u_j| + (1-\eta'(\lambda_k(u_j-u_k))) \, |\nabla u_k|} 
\end{align} 
at each point in $\hat{\Sigma} \cap W$, where $|\nabla u_j|$ and $|\nabla u_k|$ are computed with respect to the metric $g$. 

Differentiating the identity (\ref{formula.for.alpha}) gives $|\nabla^{\hat{\Sigma}} \cos(2\alpha)| \leq C$ at $x_0$, where $C$ is independent of $\gamma$ and $\lambda_0$. Since $\alpha \in (0,\frac{\pi}{2})$ is bounded away from $0$ and $\frac{\pi}{2}$ at the point $x_0$, it follows that $|\nabla^{\hat{\Sigma}} \alpha| \leq C$ at $x_0$, where $C$ is independent of $\gamma$ and $\lambda_0$. Differentiating the identity (\ref{formula.for.theta}) gives $\nabla^{\hat{\Sigma}} \cos(2\theta\alpha) = 0$ at $x_0$. Since $\theta\alpha \in (0,\frac{\pi}{2})$, it follows that $\nabla^{\hat{\Sigma}}(\theta\alpha) = 0$ at $x_0$. This implies $|\nabla^{\hat{\Sigma}} \theta| \leq C$ at $x_0$, where $C$ is independent of $\gamma$ and $\lambda_0$. Finally, differentiating the identity (\ref{formula.for.varphi}) gives 
\[\Big | \nabla^{\hat{\Sigma}} \Big ( \frac{\tan(\varphi)}{\tan(\alpha)} \Big ) - \lambda_k \Upsilon \, (\nabla^{\hat{\Sigma}} u_j - \nabla^{\hat{\Sigma}} u_k) \Big | \leq C\] 
at $x_0$. Here, $C$ is independent of $\gamma$ and $\lambda_0$, and $\Upsilon$ is defined by 
\[\Upsilon = \frac{4\eta''(\lambda_k(u_j-u_k)) \, |\nabla u_j| \, |\nabla u_k|}{\big ( (1+\eta'(\lambda_k(u_j-u_k))) \, |\nabla u_j| + (1-\eta'(\lambda_k(u_j-u_k))) \, |\nabla u_k| \big )^2}.\] 
Note that $0 \leq \Upsilon \leq C$ at the point $x_0$, where $C$ is independent of $\gamma$ and $\lambda_0$. This implies 
\begin{equation} 
\label{gradient.of.varphi}
|\nabla^{\hat{\Sigma}} \varphi - \lambda_k \Upsilon \, \tan(\alpha) \cos^2(\varphi) \, (\nabla^{\hat{\Sigma}} u_j - \nabla^{\hat{\Sigma}} u_k)| \leq C 
\end{equation}
at the point $x_0$, where $C$ is independent of $\gamma$ and $\lambda_0$. 

By Proposition \ref{recursive.relation.for.hat.nu} (iii), the unit normal to $\hat{\Sigma}$ is given by 
\begin{equation} 
\label{formula.for.hat.nu}
\hat{\nu} = \frac{\sin(\alpha+\varphi) \, \nu_j + \sin(\alpha-\varphi) \, \nu_k}{\sin(2\alpha)} 
\end{equation}
at each point in $\hat{\Sigma} \cap W$. Differentiating the identity (\ref{formula.for.hat.nu}) gives 
\begin{equation} 
\label{estimate.for.covariant.derivative.of.hat.nu}
\Big | D_\xi \hat{\nu} - d\varphi(\xi) \, \frac{\cos(\alpha+\varphi) \, \nu_j - \cos(\alpha-\varphi) \, \nu_k}{\sin(2\alpha)} \Big | \leq C \, |\xi| 
\end{equation}
at $x_0$, where $\xi$ denotes an arbitrary vector in $T_{x_0} \hat{\Sigma}$ and $C$ is independent of $\gamma$ and $\lambda_0$. In view of the identities 
\begin{align*} 
&\langle \cos(\alpha+\varphi) \, \nu_j - \cos(\alpha-\varphi) \, \nu_k,\sin(\alpha+\varphi) \, \nu_j + \sin(\alpha-\varphi) \, \nu_k \rangle \\ 
&= \cos(\alpha+\varphi) \sin(\alpha+\varphi) - \cos(\alpha-\varphi) \sin(\alpha-\varphi) \\ 
&+ \cos(\alpha+\varphi) \sin(\alpha-\varphi) \cos(2\alpha) - \cos(\alpha-\varphi) \sin(\alpha+\varphi) \cos(2\alpha) \\ 
&= 0 
\end{align*} 
and 
\begin{align*} 
&|\cos(\alpha+\varphi) \, \nu_j - \cos(\alpha-\varphi) \, \nu_k|^2 \\ 
&= \cos^2(\alpha+\varphi) + \cos^2(\alpha-\varphi) - 2 \cos(\alpha+\varphi) \cos(\alpha-\varphi) \cos(2\alpha) \\ 
&= \sin^2(2\alpha), 
\end{align*} 
the vector 
\[\frac{\cos(\alpha+\varphi) \, \nu_j - \cos(\alpha-\varphi) \, \nu_k}{\sin(2\alpha)}\] 
lies in the tangent space $T_{x_0} \hat{\Sigma}$ and has unit length with respect to the metric $g$. The estimate (\ref{estimate.for.covariant.derivative.of.hat.nu}) implies 
\begin{equation} 
\label{lower.bound.for.H.1}
H \geq \Big \langle \nabla^{\hat{\Sigma}} \varphi,\frac{\cos(\alpha+\varphi) \, \nu_j - \cos(\alpha-\varphi) \, \nu_k}{\sin(2\alpha)} \Big \rangle - C 
\end{equation}
at $x_0$, where $C$ is independent of $\gamma$ and $\lambda_0$. Combining (\ref{gradient.of.varphi}) and (\ref{lower.bound.for.H.1}), we obtain 
\begin{align} 
\label{lower.bound.for.H.2}
H &\geq \lambda_k \Upsilon \, \tan(\alpha) \cos^2(\varphi) \notag \\ 
&\qquad \cdot \Big \langle \nabla^{\hat{\Sigma}} u_j - \nabla^{\hat{\Sigma}} u_k,\frac{\cos(\alpha+\varphi) \, \nu_j - \cos(\alpha-\varphi) \, \nu_k}{\sin(2\alpha)} \Big \rangle \\ 
&- C \notag 
\end{align}
at $x_0$, where $C$ is independent of $\gamma$ and $\lambda_0$. Using the identity (\ref{formula.for.hat.nu}), we obtain 
\[\nu_j = \sin(\alpha-\varphi) \, \frac{\cos(\alpha+\varphi) \, \nu_j - \cos(\alpha-\varphi) \, \nu_k}{\sin(2\alpha)} + \cos(\alpha-\varphi) \, \hat{\nu}\] 
and 
\[\nu_k = -\sin(\alpha+\varphi) \, \frac{\cos(\alpha+\varphi) \, \nu_j - \cos(\alpha-\varphi) \, \nu_k}{\sin(2\alpha)} + \cos(\alpha+\varphi) \, \hat{\nu}.\] 
This implies 
\begin{align*} 
&\nabla u_j - \nabla u_k \\ 
&= |\nabla u_j| \, \nu_j - |\nabla u_k| \, \nu_k \\ 
&= (\sin(\alpha-\varphi) \, |\nabla u_j| + \sin(\alpha+\varphi) \, |\nabla u_k|) \, \frac{\cos(\alpha+\varphi) \, \nu_j - \cos(\alpha-\varphi) \, \nu_k}{\sin(2\alpha)} \\ 
&+ (\cos(\alpha-\varphi) \, |\nabla u_j| - \cos(\alpha+\varphi) \, |\nabla u_k|) \, \hat{\nu}. 
\end{align*} 
Projecting to $T_{x_0} \hat{\Sigma}$ gives 
\begin{align*} 
&\nabla^{\hat{\Sigma}} u_j - \nabla^{\hat{\Sigma}} u_k \\ 
&= (\sin(\alpha-\varphi) \, |\nabla u_j| + \sin(\alpha+\varphi) \, |\nabla u_k|) \, \frac{\cos(\alpha+\varphi) \, \nu_j - \cos(\alpha-\varphi) \, \nu_k}{\sin(2\alpha)}. 
\end{align*}
Since $\sin(\alpha-\varphi) \, |\nabla u_j| + \sin(\alpha+\varphi) \, |\nabla u_k| \geq 0$, it follows that 
\[\Big \langle \nabla^{\hat{\Sigma}} u_j - \nabla^{\hat{\Sigma}} u_k,\frac{\cos(\alpha+\varphi) \, \nu_j - \cos(\alpha-\varphi) \, \nu_k}{\sin(2\alpha)} \Big \rangle = |\nabla^{\hat{\Sigma}} u_j - \nabla^{\hat{\Sigma}} u_k|.\] 
Substituting this identity into (\ref{lower.bound.for.H.2}) gives 
\begin{equation} 
\label{lower.bound.for.H.3}
H \geq \lambda_k \Upsilon \, \tan(\alpha) \cos^2(\varphi) \, |\nabla^{\hat{\Sigma}} u_j - \nabla^{\hat{\Sigma}} u_k| - C 
\end{equation}
at $x_0$, where $C$ is independent of $\gamma$ and $\lambda_0$.

On the other hand, it follows from Proposition \ref{recursive.relation.for.hat.N} (vii) that 
\begin{equation} 
\label{formula.for.hat.N}
\hat{N} = \frac{\sin(\theta(\alpha+\varphi)) \, N_j + \sin(\theta(\alpha-\varphi)) \, N_k}{\sin(2\theta\alpha)} 
\end{equation} 
at each point in $\hat{\Sigma} \cap W$. Differentiating the identity (\ref{formula.for.hat.N}) gives 
\begin{equation} 
\label{estimate.for.derivative.of.hat.N}
\Big | d\hat{N}(\xi) - \theta \, d\varphi(\xi) \, \frac{\cos(\theta(\alpha+\varphi)) \, N_j - \cos(\theta(\alpha-\varphi)) \, N_k}{\sin(2\theta\alpha)} \Big | \leq C \, |\xi| 
\end{equation}
at $x_0$, where $\xi$ denotes an arbitrary vector in $T_{x_0} \hat{\Sigma}$ and $C$ is independent of $\gamma$ and $\lambda_0$. In view of the identity 
\begin{align*} 
&|\cos(\theta(\alpha+\varphi)) \, N_j - \cos(\theta(\alpha-\varphi)) \, N_k|^2 \\ 
&= \cos^2(\theta(\alpha+\varphi)) + \cos^2(\theta(\alpha-\varphi)) - 2 \cos(\theta(\alpha+\varphi)) \cos(\theta(\alpha-\varphi)) \cos(2\theta\alpha) \\ 
&= \sin^2(2\theta\alpha), 
\end{align*} 
the vector 
\[\frac{\cos(\theta(\alpha+\varphi)) \, N_j - \cos(\theta(\alpha-\varphi)) \, N_k}{\sin(2\theta\alpha)} \in \mathbb{R}^n\] 
has unit length with respect to the Euclidean metric. Therefore, the estimate (\ref{estimate.for.derivative.of.hat.N}) implies 
\begin{equation} 
\label{formula.for.trace.norm.of.differential.of.hat.N.1}
\|d\hat{N}\|_{\text{\rm tr}} \leq \theta \, |\nabla^{\hat{\Sigma}} \varphi| + C 
\end{equation}
at $x_0$, where $C$ is independent of $\gamma$ and $\lambda_0$. Combining (\ref{gradient.of.varphi}) and (\ref{formula.for.trace.norm.of.differential.of.hat.N.1}), we obtain 
\begin{equation} 
\label{formula.for.trace.norm.of.differential.of.hat.N.2}
\|d\hat{N}\|_{\text{\rm tr}} \leq \theta \lambda_k \Upsilon \, \tan(\alpha) \cos^2(\varphi) \, |\nabla^{\hat{\Sigma}} u_j - \nabla^{\hat{\Sigma}} u_k| + C 
\end{equation}
at $x_0$, where $C$ is independent of $\gamma$ and $\lambda_0$. Combining (\ref{lower.bound.for.H.3}) and (\ref{formula.for.trace.norm.of.differential.of.hat.N.2}), we conclude that 
\begin{align*} 
H - \|d\hat{N}\|_{\text{\rm tr}} 
&\geq (1-\theta) \, \lambda_k \Upsilon \, \tan(\alpha) \cos^2(\varphi) \, |\nabla^{\hat{\Sigma}} u_j - \nabla^{\hat{\Sigma}} u_k| - C \\ 
&\geq -C\lambda_k \max \{\theta-1,0\} - C 
\end{align*}
at $x_0$, where $C$ is independent of $\gamma$ and $\lambda_0$. This completes the proof of Proposition \ref{pointwise.bound.on.E_jk}. \\

\begin{corollary}
\label{pointwise.bound.on.E_jk.2}
Let $0 \leq j < k \leq q$. If $\lambda_0$ is sufficiently large (depending on $\gamma$), then $H - \|d\hat{N}\|_{\text{\rm tr}} \geq -C\gamma\lambda_k$ on the set $E_{j,k}$. Here, $C$ is independent of $\gamma$ and $\lambda_0$. 
\end{corollary}

\textbf{Proof.} 
We apply Lemma \ref{approximate.angle.inequality.2} with $\varepsilon = \gamma$. Consequently, if $\lambda_0$ is sufficiently large (depending on $\gamma$), then $\langle \nu_j,\nu_k \rangle - \langle N_j,N_k \rangle \leq \gamma$ at each point in $\Omega \cap \{-2\lambda_0^{-1} \leq u_j \leq 0\} \cap \{-2\lambda_0^{-1} \leq u_k \leq 0\}$. Since $E_{j,k} \subset \Omega \cap \{-2\lambda_k^{-1} \leq u_j \leq 0\} \cap \{-2\lambda_k^{-1} \leq u_k \leq 0\}$, it follows that $\langle \nu_j,\nu_k \rangle - \langle N_j,N_k \rangle \leq \gamma$ at each point in $E_{j,k}$. The assertion follows now from Proposition \ref{pointwise.bound.on.E_jk}. This completes the proof of Corollary \ref{pointwise.bound.on.E_jk.2}. \\

\begin{corollary}
\label{pointwise.bound.on.E_jk.3}
We can find a constant $M>2$ (independent of $\gamma$ and $\lambda_0$) with the following property. Let $0 \leq j < k \leq q$. If $\lambda_0$ is sufficiently large, then $H - \|d\hat{N}\|_{\text{\rm tr}} \geq -C$ on the set $E_{j,k} \setminus \bigcup_{i \in \{0,1,\hdots,q\} \setminus \{j,k\}} \{-M\lambda_k^{-1} \leq u_i \leq 0\}$. Here, $C$ is independent of $\gamma$ and $\lambda_0$. 
\end{corollary}

\textbf{Proof.} 
We distinguish two cases: 

\textit{Case 1:} Suppose that $\Omega \cap \{u_j=0\} \cap \{u_k=0\} = \emptyset$. Then $\Omega \cap \{-2\lambda_0^{-1} \leq u_k \leq 0\} \cap \{-2\lambda_0^{-1} \leq u_j \leq 0\} = \emptyset$ if $\lambda_0$ is sufficiently large. Consequently, $E_{j,k} = \emptyset$ if $\lambda_0$ is sufficiently large. Thus, the assertion is trivially true in this case.

\textit{Case 2:} Suppose that $\Omega \cap \{u_j=0\} \cap \{u_k=0\} \neq \emptyset$. It follows from Assumption \ref{no.redundant.inequalities} that the hyperplanes $\{u_j=0\}$ and $\{u_k=0\}$ must intersect transversally. Let us consider a point $x \in E_{j,k}$. Then $-2\lambda_k^{-1} \leq u_j(x) \leq 0$ and $-2\lambda_k^{-1} \leq u_k(x) \leq 0$. By transversality, we can find a point $y$ such that $u_j(y)=u_k(y)=0$ and $d_{\text{\rm eucl}}(x,y) \leq M\lambda_k^{-1}$, where $M>2$ is independent of $\gamma$ and $\lambda_0$. 

Suppose now that $x \in E_{j,k} \setminus \bigcup_{i \in \{0,1,\hdots,q\} \setminus \{j,k\}} \{-M\lambda_k^{-1} \leq u_i \leq 0\}$. Then $u_i(x) < -M\lambda_k^{-1}$ for all $i \in \{0,1,\hdots,q\} \setminus \{j,k\}$. Since $|u_i(x)-u_i(y)| \leq d_{\text{\rm eucl}}(x,y) \leq M\lambda_k^{-1}$ for all $i \in \{0,1,\hdots,q\}$, it follows that $u_i(y) < 0$ for all $i \in \{0,1,\hdots,q\} \setminus \{j,k\}$, and consequently $y \in \Omega$. Assumption \ref{angle.comparison} implies that $\langle \nu_j,\nu_k \rangle - \langle N_j,N_k \rangle \leq 0$ at the point $y$. Since $d_{\text{\rm eucl}}(x,y) \leq M\lambda_k^{-1}$, it follows that $\langle \nu_j,\nu_k \rangle - \langle N_j,N_k \rangle \leq C\lambda_k^{-1}$ at the point $x$, where $C$ is independent of $\gamma$ and $\lambda_0$. Using Proposition \ref{pointwise.bound.on.E_jk}, we conclude that $H - \|d\hat{N}\|_{\text{\rm tr}} \geq -C$ at the point $x$, where $C$ is independent of $\gamma$ and $\lambda_0$. This completes the proof of Corollary \ref{pointwise.bound.on.E_jk.3}. \\

\begin{proposition}
\label{pointwise.bound.on.G_ijk}
Let $0 \leq i < j < k \leq q$. If $\lambda_0$ is sufficiently large, then $H - \|d\hat{N}\|_{\text{\rm tr}} \geq -C\lambda_k \max \{\langle \hat{\nu}_j,\nu_k \rangle - \langle \hat{N}_j,N_k \rangle,0\} - C\lambda_j$ on the set $G_{i,j,k}$. Here, $C$ is independent of $\gamma$ and $\lambda_0$. 
\end{proposition}

\textbf{Proof.} 
Let us fix a point $x_0 \in G_{i,j,k}$. In the following, we will consider a small open neighborhood $W$ of $x_0$. We choose $W$ so that 
\[W \subset \bigcap_{m=k+1}^q \{\hat{u}_{m-1}-u_m > \lambda_m^{-1}\} \cap \bigcap_{m=j+1}^{k-1} \{\hat{u}_{m-1}-u_m > \lambda_m^{-1}\}.\] 
Then $\hat{u}_j = \hat{u}_{k-1}$ and $\hat{u}_k = \hat{u}_q$ at each point in $W$. This implies $\hat{\nu}_j = \hat{\nu}_{k-1}$ and $\hat{\nu}_k = \hat{\nu}_q = \hat{\nu}$ at each point in $\hat{\Sigma} \cap W$. Moreover, it follows from Proposition \ref{recursive.relation.for.hat.N} (v) that $\hat{N}_j = \hat{N}_{k-1}$ and $\hat{N}_k = \hat{N}_q = \hat{N}$ at each point in $\hat{\Sigma} \cap W$. 

We define a smooth function $\alpha: \hat{\Sigma} \cap W \to (0,\frac{\pi}{2})$ by 
\begin{equation} 
\label{formula.for.alpha.new}
\cos(2\alpha) = \langle \hat{\nu}_j,\nu_k \rangle 
\end{equation} 
at each point in $\hat{\Sigma} \cap W$. Moreover, we define a smooth function $\theta: \hat{\Sigma} \cap W \to \mathbb{R}$ so that $\theta \in (0,\frac{\pi}{2\alpha})$ at each point in $\hat{\Sigma} \cap W$ and 
\begin{equation} 
\label{formula.for.theta.new}
\cos(2\theta\alpha) = \langle \hat{N}_j,N_k \rangle 
\end{equation} 
at each point in $\hat{\Sigma} \cap W$. Finally, we define a function $\varphi: \hat{\Sigma} \cap W \to \mathbb{R}$ so that $\varphi \in [-\alpha,\alpha]$ at each point in $\hat{\Sigma} \cap W$ and 
\begin{align} 
\label{formula.for.varphi.new}
&\frac{\tan(\varphi)}{\tan(\alpha)} \notag \\ 
&= \frac{(1+\eta'(\lambda_k(\hat{u}_j-u_k))) \, |\nabla \hat{u}_j| - (1-\eta'(\lambda_k(\hat{u}_j-u_k))) \, |\nabla u_k|}{(1+\eta'(\lambda_k(\hat{u}_j-u_k))) \, |\nabla \hat{u}_j| + (1-\eta'(\lambda_k(\hat{u}_j-u_k))) \, |\nabla u_k|} 
\end{align} 
at each point in $\hat{\Sigma} \cap W$, where $|\nabla \hat{u}_j|$ and $|\nabla u_k|$ are computed with respect to the metric $g$. 

Let us differentiate the identity (\ref{formula.for.alpha.new}). Using Lemma \ref{C2.bound.for.hat.u} and Lemma \ref{lower.bound.for.norm.of.gradient.of.hat.u_k}, we obtain $|\nabla^{\hat{\Sigma}} \cos(2\alpha)| \leq C \lambda_j$ at $x_0$, where $C$ is independent of $\gamma$ and $\lambda_0$. Since $\alpha \in (0,\frac{\pi}{2})$ is bounded away from $0$ and $\frac{\pi}{2}$ at the point $x_0$, it follows that $|\nabla^{\hat{\Sigma}} \alpha| \leq C \lambda_j$ at $x_0$, where $C$ is independent of $\gamma$ and $\lambda_0$. We next differentiate the identity (\ref{formula.for.theta.new}). Using Proposition \ref{C1.bound.for.hat.N}, we obtain $|\nabla^{\hat{\Sigma}} \cos(2\theta\alpha)| \leq C \lambda_j$ at $x_0$. Since $\theta\alpha \in (0,\frac{\pi}{2})$ is bounded away from $0$ and $\frac{\pi}{2}$ at the point $x_0$, it follows that $|\nabla^{\hat{\Sigma}}(\theta\alpha)| \leq C \lambda_j$ at $x_0$, where $C$ is independent of $\gamma$ and $\lambda_0$. This implies $|\nabla^{\hat{\Sigma}} \theta| \leq C \lambda_j$ at $x_0$, where $C$ is independent of $\gamma$ and $\lambda_0$. Finally, we differentiate the identity (\ref{formula.for.varphi.new}). Using Lemma \ref{C2.bound.for.hat.u} and Lemma \ref{lower.bound.for.norm.of.gradient.of.hat.u_k}, we obtain 
\[\Big | \nabla^{\hat{\Sigma}} \Big ( \frac{\tan(\varphi)}{\tan(\alpha)} \Big ) - \lambda_k \Upsilon \, (\nabla^{\hat{\Sigma}} \hat{u}_j - \nabla^{\hat{\Sigma}} u_k) \Big | \leq C \lambda_j\] 
at $x_0$. Here, $C$ is independent of $\gamma$ and $\lambda_0$, and $\Upsilon$ is defined by 
\[\Upsilon = \frac{4\eta''(\lambda_k(\hat{u}_j-u_k)) \, |\nabla \hat{u}_j| \, |\nabla u_k|}{\big ( (1+\eta'(\lambda_k(\hat{u}_j-u_k))) \, |\nabla \hat{u}_j| + (1-\eta'(\lambda_k(\hat{u}_j-u_k))) \, |\nabla u_k| \big )^2}.\] 
Note that $0 \leq \Upsilon \leq C$ at the point $x_0$, where $C$ is independent of $\gamma$ and $\lambda_0$. This implies 
\begin{equation} 
\label{gradient.of.varphi.new}
|\nabla^{\hat{\Sigma}} \varphi - \lambda_k \Upsilon \, \tan(\alpha) \cos^2(\varphi) \, (\nabla^{\hat{\Sigma}} \hat{u}_j - \nabla^{\hat{\Sigma}} u_k)| \leq C \lambda_j
\end{equation}
at the point $x_0$, where $C$ is independent of $\gamma$ and $\lambda_0$. 

By Proposition \ref{recursive.relation.for.hat.nu} (iii), the unit normal to $\hat{\Sigma}$ is given by 
\begin{equation} 
\label{formula.for.hat.nu.new}
\hat{\nu} = \frac{\sin(\alpha+\varphi) \, \hat{\nu}_j + \sin(\alpha-\varphi) \, \nu_k}{\sin(2\alpha)} 
\end{equation}
at each point in $\hat{\Sigma} \cap W$. Differentiating the identity (\ref{formula.for.hat.nu.new}) gives 
\begin{equation} 
\label{estimate.for.covariant.derivative.of.hat.nu.new}
\Big | D_\xi \hat{\nu} - d\varphi(\xi) \, \frac{\cos(\alpha+\varphi) \, \hat{\nu}_j - \cos(\alpha-\varphi) \, \nu_k}{\sin(2\alpha)} \Big | \leq C \lambda_j \, |\xi| 
\end{equation} 
at $x_0$, where $\xi$ denotes an arbitrary vector in $T_{x_0} \hat{\Sigma}$ and $C$ is independent of $\gamma$ and $\lambda_0$. The vector 
\[\frac{\cos(\alpha+\varphi) \, \hat{\nu}_j - \cos(\alpha-\varphi) \, \nu_k}{\sin(2\alpha)}\] 
lies in the tangent space $T_{x_0} \hat{\Sigma}$ and has unit length with respect to the metric $g$. The estimate (\ref{estimate.for.covariant.derivative.of.hat.nu.new}) implies 
\begin{equation} 
\label{lower.bound.for.H.1.new}
H \geq \Big \langle \nabla^{\hat{\Sigma}} \varphi,\frac{\cos(\alpha+\varphi) \, \hat{\nu}_j - \cos(\alpha-\varphi) \, \nu_k}{\sin(2\alpha)} \Big \rangle - C \lambda_j 
\end{equation}
at $x_0$, where $C$ is independent of $\gamma$ and $\lambda_0$. Combining (\ref{gradient.of.varphi.new}) and (\ref{lower.bound.for.H.1.new}), we obtain 
\begin{align} 
\label{lower.bound.for.H.2.new}
H &\geq \lambda_k \Upsilon \, \tan(\alpha) \cos^2(\varphi) \notag \\ 
&\qquad \cdot \Big \langle \nabla^{\hat{\Sigma}} \hat{u}_j - \nabla^{\hat{\Sigma}} u_k,\frac{\cos(\alpha+\varphi) \, \hat{\nu}_j - \cos(\alpha-\varphi) \, \nu_k}{\sin(2\alpha)} \Big \rangle \\ 
&- C \lambda_j \notag 
\end{align}
at $x_0$, where $C$ is independent of $\gamma$ and $\lambda_0$. Using the identity (\ref{formula.for.hat.nu.new}), we obtain 
\[\hat{\nu}_j = \sin(\alpha-\varphi) \, \frac{\cos(\alpha+\varphi) \, \hat{\nu}_j - \cos(\alpha-\varphi) \, \nu_k}{\sin(2\alpha)} + \cos(\alpha-\varphi) \, \hat{\nu}\] 
and 
\[\nu_k = -\sin(\alpha+\varphi) \, \frac{\cos(\alpha+\varphi) \, \hat{\nu}_j - \cos(\alpha-\varphi) \, \nu_k}{\sin(2\alpha)} + \cos(\alpha+\varphi) \, \hat{\nu}.\] 
This implies 
\begin{align*} 
&\nabla \hat{u}_j - \nabla u_k \\ 
&= |\nabla \hat{u}_j| \, \hat{\nu}_j - |\nabla u_k| \, \nu_k \\ 
&= (\sin(\alpha-\varphi) \, |\nabla \hat{u}_j| + \sin(\alpha+\varphi) \, |\nabla u_k|) \, \frac{\cos(\alpha+\varphi) \, \hat{\nu}_j - \cos(\alpha-\varphi) \, \nu_k}{\sin(2\alpha)} \\ 
&+ (\cos(\alpha-\varphi) \, |\nabla \hat{u}_j| - \cos(\alpha+\varphi) \, |\nabla u_k|) \, \hat{\nu}. 
\end{align*} 
Projecting to $T_{x_0} \hat{\Sigma}$ gives 
\begin{align*} 
&\nabla^{\hat{\Sigma}} \hat{u}_j - \nabla^{\hat{\Sigma}} u_k \\ 
&= (\sin(\alpha-\varphi) \, |\nabla \hat{u}_j| + \sin(\alpha+\varphi) \, |\nabla u_k|) \, \frac{\cos(\alpha+\varphi) \, \hat{\nu}_j - \cos(\alpha-\varphi) \, \nu_k}{\sin(2\alpha)}. 
\end{align*}
Since $\sin(\alpha-\varphi) \, |\nabla \hat{u}_j| + \sin(\alpha+\varphi) \, |\nabla u_k| \geq 0$, it follows that 
\[\Big \langle \nabla^{\hat{\Sigma}} \hat{u}_j - \nabla^{\hat{\Sigma}} u_k,\frac{\cos(\alpha+\varphi) \, \hat{\nu}_j - \cos(\alpha-\varphi) \, \nu_k}{\sin(2\alpha)} \Big \rangle = |\nabla^{\hat{\Sigma}} \hat{u}_j - \nabla^{\hat{\Sigma}} u_k|.\] 
Substituting this identity into (\ref{lower.bound.for.H.2.new}) gives 
\begin{equation} 
\label{lower.bound.for.H.3.new}
H \geq \lambda_k \Upsilon \, \tan(\alpha) \cos^2(\varphi) \, |\nabla^{\hat{\Sigma}} \hat{u}_j - \nabla^{\hat{\Sigma}} u_k| - C \lambda_j 
\end{equation}
at $x_0$, where $C$ is independent of $\gamma$ and $\lambda_0$.

On the other hand, it follows from Proposition \ref{recursive.relation.for.hat.N} (vii) that 
\begin{equation} 
\label{formula.for.hat.N.new}
\hat{N} = \frac{\sin(\theta(\alpha+\varphi)) \, \hat{N}_j + \sin(\theta(\alpha-\varphi)) \, N_k}{\sin(2\theta\alpha)} 
\end{equation} 
at each point in $\hat{\Sigma} \cap W$. Differentiating the identity (\ref{formula.for.hat.N.new}) gives 
\begin{equation} 
\label{estimate.for.derivative.of.hat.N.new}
\Big | d\hat{N}(\xi) - \theta \, d\varphi(\xi) \, \frac{\cos(\theta(\alpha+\varphi)) \, \hat{N}_j - \cos(\theta(\alpha-\varphi)) \, N_k}{\sin(2\theta\alpha)} \Big | \leq C \lambda_j \, |\xi| 
\end{equation}
at $x_0$, where $\xi$ denotes an arbitrary vector in $T_{x_0} \hat{\Sigma}$ and $C$ is independent of $\gamma$ and $\lambda_0$. The vector 
\[\frac{\cos(\theta(\alpha+\varphi)) \, \hat{N}_j - \cos(\theta(\alpha-\varphi)) \, N_k}{\sin(2\theta\alpha)} \in \mathbb{R}^n\] 
has unit length with respect to the Euclidean metric. Therefore, the estimate (\ref{estimate.for.derivative.of.hat.N.new}) implies 
\begin{equation} 
\label{formula.for.trace.norm.of.dhatN.1.new}
\|d\hat{N}\|_{\text{\rm tr}} \leq \theta \, |\nabla^{\hat{\Sigma}} \varphi| + C \lambda_j 
\end{equation}
at $x_0$, where $C$ is independent of $\gamma$ and $\lambda_0$. Combining (\ref{gradient.of.varphi.new}) and (\ref{formula.for.trace.norm.of.dhatN.1.new}), we obtain 
\begin{equation} 
\label{formula.for.trace.norm.of.dhatN.2.new}
\|d\hat{N}\|_{\text{\rm tr}} \leq \theta \lambda_k \Upsilon \, \tan(\alpha) \cos^2(\varphi) \, |\nabla^{\hat{\Sigma}} \hat{u}_j - \nabla^{\hat{\Sigma}} u_k| + C \lambda_j 
\end{equation}
at $x_0$, where $C$ is independent of $\gamma$ and $\lambda_0$. Combining (\ref{lower.bound.for.H.3.new}) and (\ref{formula.for.trace.norm.of.dhatN.2.new}), we conclude that 
\begin{align*} 
H - \|d\hat{N}\|_{\text{\rm tr}} 
&\geq (1-\theta) \, \lambda_k \Upsilon \, \tan(\alpha) \cos^2(\varphi) \, |\nabla^{\hat{\Sigma}} \hat{u}_j - \nabla^{\hat{\Sigma}} u_k| - C \lambda_j \\ 
&\geq -C\lambda_k \max \{\theta-1,0\} - C \lambda_j 
\end{align*}
at $x_0$, where $C$ is independent of $\gamma$ and $\lambda_0$. This completes the proof of Proposition \ref{pointwise.bound.on.G_ijk}. \\

\begin{corollary}
\label{pointwise.bound.on.G_ijk.2}
Let $0 \leq i < j < k \leq q$. If $\lambda_0$ is sufficiently large (depending on $\gamma$), then $H - \|d\hat{N}\|_{\text{\rm tr}} \geq -C\gamma\lambda_k$ on the set $G_{i,j,k}$. Here, $C$ is independent of $\gamma$ and $\lambda_0$. 
\end{corollary}

\textbf{Proof.} 
We apply Proposition \ref{approximate.angle.inequality.3} with $\varepsilon = \gamma$. Consequently, if $\lambda_0$ is sufficiently large (depending on $\gamma$), then $\langle \hat{\nu}_j,\nu_k \rangle - \langle \hat{N}_j,N_k \rangle \leq \gamma$ at each point in $\Omega \cap \{-2\lambda_0^{-1} \leq \hat{u}_j \leq 0\} \cap \{-2\lambda_0^{-1} \leq u_k \leq 0\}$. Since $G_{i,j,k} \subset \Omega \cap \{-2\lambda_k^{-1} \leq \hat{u}_j \leq 0\} \cap \{-2\lambda_k^{-1} \leq u_k \leq 0\}$, it follows that $\langle \hat{\nu}_j,\nu_k \rangle - \langle \hat{N}_j,N_k \rangle \leq \gamma$ at each point in $G_{i,j,k}$. The assertion follows now from Proposition \ref{pointwise.bound.on.G_ijk}. This completes the proof of Corollary \ref{pointwise.bound.on.G_ijk.2}. \\

In the remainder of this section, we will estimate a suitable Morrey norm of the function $\max \{\|d\hat{N}\|_{\text{\rm tr}} - H,0\}$ on $\hat{\Sigma}$. To that end, we will combine Proposition \ref{pointwise.bound.on.F_k}, Corollary \ref{pointwise.bound.on.E_jk.2}, Corollary \ref{pointwise.bound.on.E_jk.3}, and Corollary \ref{pointwise.bound.on.G_ijk.2} with the area estimates in Lemma \ref{area.1} and Lemma \ref{area.2}. 

\begin{proposition} 
\label{Morrey.norm.E_jk}
Fix an exponent $\sigma \in (1,\frac{q}{q-1})$. Suppose that $\gamma \in (0,M^{-1})$, where $M>2$ denotes the constant in Corollary \ref{pointwise.bound.on.E_jk.3}. If $\lambda_0$ is sufficiently large (depending on $\gamma$), then 
\[r^{\sigma+1-n} \int_{E_{j,k} \cap B_r(p)} (\max \{\|d\hat{N}\|_{\text{\rm tr}} - H,0\})^\sigma \leq C \, \gamma^{\sigma-q(\sigma-1)}\] 
for all $0 \leq j < k \leq q$ and all $0 < r \leq 1$. Here, $C$ is independent of $\gamma$ and $\lambda_0$. 
\end{proposition}

\textbf{Proof.} 
If $\lambda_0$ is sufficiently large (depending on $\gamma$), then Corollary \ref{pointwise.bound.on.E_jk.2} implies that $\max \{\|d\hat{N}\|_{\text{\rm tr}} - H,0\} \leq C\gamma\lambda_k$ on the set $E_{j,k}$, where $C$ is independent of $\gamma$ and $\lambda_0$. Moreover, Corollary \ref{pointwise.bound.on.E_jk.3} implies that $\max \{\|d\hat{N}\|_{\text{\rm tr}} - H,0\} \leq C$ on the set $E_{j,k} \setminus \bigcup_{i \in \{0,1,\hdots,q\} \setminus \{j,k\}} \{-M\lambda_k^{-1} \leq u_i \leq 0\}$, where, again, $C$ is independent of $\gamma$ and $\lambda_0$. We distinguish two cases: 

\textit{Case 1:} Suppose that $0 < r \leq \lambda_0^{-1}$. Since 
\[E_{j,k} \subset \Omega \cap \{\hat{u}_k=0\} \cap \{-2\lambda_k^{-1} \leq \hat{u}_{k-1} \leq 0\} \cap \{-2\lambda_k^{-1} \leq u_k \leq 0\},\] 
Lemma \ref{area.1} implies that $|E_{j,k} \cap B_r(p)| \leq C \lambda_k^{-1} r^{n-2}$, where $C$ is independent of $\gamma$ and $\lambda_0$. Consequently, 
\begin{align*} 
&r^{\sigma+1-n} \int_{E_{j,k} \cap B_r(p)} (\max \{\|d\hat{N}\|_{\text{\rm tr}} - H,0\})^\sigma \\ 
&\leq C \, r^{\sigma+1-n} \, (\gamma\lambda_k)^\sigma \, |E_{j,k} \cap B_r(p)| \\ 
&\leq C \, \gamma^\sigma \, (\lambda_k r)^{\sigma-1} \\ 
&\leq C \, \gamma^{\sigma-q(\sigma-1)}, 
\end{align*}
where $C$ is independent of $\gamma$ and $\lambda_0$.

\textit{Case 2:} Suppose that $\lambda_0^{-1} \leq r \leq 1$. As above, $|E_{j,k} \cap B_r(p)| \leq C \lambda_k^{-1} r^{n-2}$ by Lemma \ref{area.1}. Moreover, since $\gamma \in (0,M^{-1})$, we have $M\lambda_k^{-1} \leq (\gamma\lambda_k)^{-1} \leq \lambda_0^{-1}$. This implies 
\begin{align*} 
&E_{j,k} \cap \{-M\lambda_k^{-1} \leq u_i \leq 0\} \\ 
&\subset \Omega \cap \{\hat{u}_k=0\} \cap \{-2\lambda_k^{-1} \leq \hat{u}_{k-1} \leq 0\} \cap \{-2\lambda_k^{-1} \leq u_k \leq 0\} \\ 
&\cap \{-4\lambda_0^{-1} \leq u_j \leq 0\} \cap \{-6\lambda_0^{-1} \leq u_i \leq 0\} 
\end{align*}
for each $i \in \{0,1,\hdots,q\} \setminus \{j,k\}$. Using Lemma \ref{area.2}, we obtain 
\[|E_{j,k} \cap \{-M\lambda_k^{-1} \leq u_i \leq 0\} \cap B_r(p)| \leq C \lambda_k^{-1} \lambda_0^{-1} r^{n-3}\] 
for each $i \in \{0,1,\hdots,q\} \setminus \{j,k\}$, where $C$ is independent of $\gamma$ and $\lambda_0$. Consequently, 
\begin{align*} 
&r^{\sigma+1-n} \int_{E_{j,k} \cap B_r(p)} (\max \{\|d\hat{N}\|_{\text{\rm tr}} - H,0\})^\sigma \\ 
&\leq C \, r^{\sigma+1-n} \, (\gamma\lambda_k)^\sigma \,  \sum_{i \in \{0,1,\hdots,q\} \setminus \{j,k\}} |E_{j,k} \cap \{-M\lambda_k^{-1} \leq u_i \leq 0\} \cap B_r(p)| \\ 
&+ C \, r^{\sigma+1-n} \, |E_{j,k} \cap B_r(p)| \\ 
&\leq C \, \gamma^\sigma \, (\lambda_k \lambda_0^{-1})^{\sigma-1} \, (\lambda_0 r)^{\sigma-2} + C \lambda_k^{-1} \, r^{\sigma-1} \\ 
&\leq C \, \gamma^{\sigma-q(\sigma-1)} + C \lambda_k^{-1}, 
\end{align*}
where $C$ is independent of $\gamma$ and $\lambda_0$. This completes the proof of Proposition \ref{Morrey.norm.E_jk}. \\

\begin{proposition}
\label{Morrey.norm.G_ijk}
Fix an exponent $\sigma \in (1,\frac{q}{q-1})$. If $\lambda_0$ is sufficiently large (depending on $\gamma$), then 
\[r^{\sigma+1-n} \int_{G_{i,j,k} \cap B_r(p)} (\max \{\|d\hat{N}\|_{\text{\rm tr}} - H,0\})^\sigma \leq C \, \gamma^{\sigma-q(\sigma-1)}\] 
for all $0 \leq i < j < k \leq q$ and all $0 < r \leq 1$. Here, $C$ is independent of $\gamma$ and $\lambda_0$. 
\end{proposition}

\textbf{Proof.} 
If $\lambda_0$ is sufficiently large (depending on $\gamma$), then Corollary \ref{pointwise.bound.on.G_ijk.2} implies that $\max \{\|d\hat{N}\|_{\text{\rm tr}} - H,0\} \leq C\gamma\lambda_k$ on the set $G_{i,j,k}$, where $C$ is independent of $\gamma$ and $\lambda_0$. We distinguish two cases: 

\textit{Case 1:} Suppose that $0 < r \leq \lambda_0^{-1}$. Since 
\[G_{i,j,k} \subset \Omega \cap \{\hat{u}_k=0\} \cap \{-2\lambda_k^{-1} \leq \hat{u}_{k-1} \leq 0\} \cap \{-2\lambda_k^{-1} \leq u_k \leq 0\},\] 
Lemma \ref{area.1} implies that $|G_{i,j,k} \cap B_r(p)| \leq C \lambda_k^{-1} r^{n-2}$, where $C$ is independent of $\gamma$ and $\lambda_0$. Consequently, 
\begin{align*} 
&r^{\sigma+1-n} \int_{G_{i,j,k} \cap B_r(p)} (\max \{\|d\hat{N}\|_{\text{\rm tr}} - H,0\})^\sigma \\ 
&\leq C \, r^{\sigma+1-n} \, (\gamma \lambda_k)^\sigma \, |G_{i,j,k} \cap B_r(p)| \\ 
&\leq C \, \gamma^\sigma \, (\lambda_k r)^{\sigma-1} \\ 
&\leq C \, \gamma^{\sigma-q(\sigma-1)}, 
\end{align*}
where $C$ is independent of $\gamma$ and $\lambda_0$.

\textit{Case 2:} Suppose that $\lambda_0^{-1} \leq r \leq 1$. Since 
\begin{align*} 
G_{i,j,k} 
&\subset \Omega \cap \{\hat{u}_k=0\} \cap \{-2\lambda_k^{-1} \leq \hat{u}_{k-1} \leq 0\} \cap \{-2\lambda_k^{-1} \leq u_k \leq 0\} \\ 
&\cap \{-4\lambda_0^{-1} \leq u_j \leq 0\} \cap \{-6\lambda_0^{-1} \leq u_i \leq 0\}, 
\end{align*} 
Lemma \ref{area.2} implies that $|G_{i,j,k} \cap B_r(p)| \leq C \lambda_k^{-1} \lambda_0^{-1} r^{n-3}$, where $C$ is independent of $\gamma$ and $\lambda_0$. Consequently, 
\begin{align*} 
&r^{\sigma+1-n} \int_{G_{i,j,k} \cap B_r(p)} (\max \{\|d\hat{N}\|_{\text{\rm tr}} - H,0\})^\sigma \\ 
&\leq C \, r^{\sigma+1-n} \, (\gamma \lambda_k)^\sigma \, |G_{i,j,k} \cap B_r(p)| \\ 
&\leq C \, \gamma^\sigma \, (\lambda_k \lambda_0^{-1})^{\sigma-1} \, (\lambda_0 r)^{\sigma-2} \\ 
&\leq C \, \gamma^{\sigma-q(\sigma-1)}, 
\end{align*}
where $C$ is independent of $\gamma$ and $\lambda_0$. This completes the proof of Proposition \ref{Morrey.norm.G_ijk}. \\

\begin{corollary} 
\label{Morrey.norm}
Fix an exponent $\sigma \in (1,\frac{q}{q-1})$. Suppose that $\gamma \in (0,M^{-1})$, where $M>2$ denotes the constant in Corollary \ref{pointwise.bound.on.E_jk.3}. If $\lambda_0$ is sufficiently large (depending on $\gamma$), then 
\[r^{\sigma+1-n} \int_{\hat{\Sigma} \cap B_r(p)} (\max \{\|d\hat{N}\|_{\text{\rm tr}} - H,0\})^\sigma \leq C \, \gamma^{\sigma-q(\sigma-1)}\] 
for all $0 < r \leq 1$. Here, $C$ is independent of $\gamma$ and $\lambda_0$. 
\end{corollary}

\textbf{Proof.} 
This follows by combining Proposition \ref{decomposition.into.subsets}, Proposition \ref{pointwise.bound.on.F_k}, Proposition \ref{Morrey.norm.E_jk}, and Proposition \ref{Morrey.norm.G_ijk}. \\

\section{Proof of Theorem \ref{main.thm}} 

\label{proof.of.main.thm}

In this section, we give the proof of Theorem \ref{main.thm}. It suffices to consider the odd-dimensional case. (The even-dimensional case can be reduced to the odd-dimensional case by passing to the Cartesian product $\Omega \times [-1,1] \subset \mathbb{R}^{n+1}$.) Suppose that $n \geq 3$ is an odd integer, and $\Omega = \bigcap_{m=0}^q \{u_m \leq 0\}$ is a compact, convex polytope in $\mathbb{R}^n$ with non-empty interior. Let $g$ be a Riemannian metric on $\mathbb{R}^n$. We assume that the assumptions of Theorem \ref{main.thm} are satisfied. In other words, Assumptions \ref{no.redundant.inequalities}, \ref{gradient.of.u_k}, \ref{angles.bounded.by.pi/2}, \ref{angle.comparison}, and \ref{faces.are.mean.convex} are satisfied, and $g$ has nonnegative scalar curvature at each point in $\Omega$.

Let us fix an exponent $\sigma \in (1,\frac{q}{q-1})$. The results in the previous sections imply that we can find a sequence of domains $\hat{\Omega}^{(l)}$ and a sequence of smooth maps $\hat{N}^{(l)}: \partial \hat{\Omega}^{(l)} \to S^{n-1}$ with the following properties: 
\begin{itemize}
\item For each $l$, $\hat{\Omega}^{(l)}$ is compact, convex domain in $\mathbb{R}^n$ with smooth boundary $\partial \hat{\Omega}^{(l)} = \hat{\Sigma}^{(l)}$.
\item For each $l$, $\bigcap_{m=0}^q \{u_m \leq -l^{-1}\} \subset \hat{\Omega}^{(l)} \subset \Omega$.
\item For each $l$, the map $\hat{N}^{(l)}: \hat{\Sigma}^{(l)} \to S^{n-1}$ is homotopic to the Gauss map of $\hat{\Sigma}^{(l)}$.
\item We have 
\begin{equation} 
\label{Morrey.norm.converging.to.0}
\sup_{p \in \mathbb{R}^n} \sup_{0 < r \leq 1} r^{\sigma+1-n} \int_{\hat{\Sigma}^{(l)} \cap B_r(p)} (\max \{\|d\hat{N}^{(l)}\|_{\text{\rm tr}} - H,0\})^\sigma \to 0 
\end{equation}
as $l \to \infty$.
\end{itemize}
Let us fix a Euclidean ball $U$ with the property that the closure of $U$ is contained in the interior of $\Omega$. Note that $U \subset \hat{\Omega}^{(l)}$ if $l$ is sufficiently large. In the following, we will always assume that $l$ is chosen sufficiently large so that $U \subset \hat{\Omega}^{(l)}$.

\begin{proposition}
\label{L2}
There exists a uniform constant $C$ (independent of $l$) such that
\[\int_{\hat{\Omega}^{(l)}} F^2 \, d\text{\rm vol}_g \leq C \int_{\hat{\Omega}^{(l)}} |\nabla F|^2 \, d\text{\rm vol}_g + C \int_U F^2 \, d\text{\rm vol}_g\] 
for every smooth function $F: \hat{\Omega}^{(l)} \to \mathbb{R}$. 
\end{proposition} 

\textbf{Proof.} 
The hypersurface $\hat{\Sigma}^{(l)} = \partial \hat{\Omega}^{(l)}$ can be written as a radial graph with bounded slope. From this, it is easy to see that $\hat{\Omega}^{(l)}$ is bi-Lipschitz equivalent to the Euclidean unit ball, with constants that are independent of $l$. The assertion follows now from the corresponding estimate on the unit ball. \\

\begin{proposition}
\label{sobolev.trace.estimate} 
There exists a uniform constant $C$ (independent of $l$) such that
\[\int_{\hat{\Sigma}^{(l)}} F^2 \, d\sigma_g \leq C \int_{\hat{\Omega}^{(l)}} |\nabla F|^2 \, d\text{\rm vol}_g + C \int_{\hat{\Omega}^{(l)}} F^2 \, d\text{\rm vol}_g\] 
for every smooth function $F: \hat{\Omega}^{(l)} \to \mathbb{R}$. 
\end{proposition} 

\textbf{Proof.} 
The hypersurface $\hat{\Sigma}^{(l)} = \partial \hat{\Omega}^{(l)}$ can be written as a radial graph with bounded slope. From this, it is easy to see that $\hat{\Omega}^{(l)}$ is bi-Lipschitz equivalent to the Euclidean unit ball, with constants that are independent of $l$. The assertion follows now from the Sobolev trace theorem on the unit ball. \\

\begin{proposition}
\label{fefferman.phong}
We have 
\[\int_{\hat{\Sigma}^{(l)}} \max \{\|d\hat{N}^{(l)}\|_{\text{\rm tr}} - H,0\} \, F^2 \, d\sigma_g \leq o(1) \int_{\hat{\Omega}^{(l)}} |\nabla F|^2 \, d\text{\rm vol}_g + o(1) \int_{\hat{\Sigma}^{(l)}} F^2 \, d\sigma_g\] 
for every smooth function $F: \hat{\Omega}^{(l)} \to \mathbb{R}$. 
\end{proposition}

\textbf{Proof.} 
The hypersurface $\hat{\Sigma}^{(l)} = \partial \hat{\Omega}^{(l)}$ can be written as a radial graph with bounded slope. From this, it is easy to see that $\hat{\Omega}^{(l)}$ is bi-Lipschitz equivalent to the Euclidean unit ball, with constants that are independent of $l$. Therefore, the assertion follows from (\ref{Morrey.norm.converging.to.0}) together with an estimate of Fefferman and Phong (see \cite{Brendle}, Appendix A, and \cite{Fefferman-Phong}). This completes the proof of Proposition \ref{fefferman.phong}. \\

As in \cite{Brendle}, we consider a boundary value problem for the Dirac operator. Let $m = 2^{[\frac{n}{2}]}$. Let $\mathcal{S}$ denote the spinor bundle over $\Omega$. Note that $\mathcal{S}$ is a complex vector bundle of rank $m$ equipped with a Hermitian inner product. We define a complex vector bundle $\mathcal{E}$ over $\Omega$ by 
\[\mathcal{E} = \underbrace{\mathcal{S} \oplus \hdots \oplus \mathcal{S}}_{\text{\rm $m$ times}}.\] 
Note that $\mathcal{E}$ has rank $m^2$.

As in \cite{Brendle}, we may use the map $\hat{N}^{(l)}: \hat{\Sigma}^{(l)} \to S^{n-1}$ to define a boundary chirality $\chi^{(l)}: \mathcal{E}|_{\hat{\Sigma}^{(l)}} \to \mathcal{E}|_{\hat{\Sigma}^{(l)}}$. 

\begin{definition}[cf. \cite{Brendle}, Definition 2.2]
Let $\{E_1,\hdots,E_n\}$ denote the standard basis of $\mathbb{R}^n$. Let $\omega_1,\hdots,\omega_n \in \text{\rm End}(\mathbb{C}^m)$ be defined as in \cite{Brendle}. We define a bundle map $\chi^{(l)}: \mathcal{E}|_{\Sigma^{(l)}} \to \mathcal{E}|_{\Sigma^{(l)}}$ by 
\[(\chi^{(l)} s)_\alpha = -\sum_{a=1}^n \sum_{\beta=1}^m \langle \hat{N}^{(l)},E_a \rangle \, \omega_{a\alpha\beta} \, \nu \cdot s_\beta\] 
for $\alpha=1,\hdots,m$.
\end{definition}

Recall that the map $\hat{N}^{(l)}: \hat{\Sigma}^{(l)} \to S^{n-1}$ is homotopic to the Gauss map of $\hat{\Sigma}^{(l)}$. By Proposition 2.15 in \cite{Brendle}, we can find an $m$-tuple of spinors $s^{(l)} = (s_1^{(l)},\hdots,s_m^{(l)})$ defined on $\hat{\Omega}^{(l)}$ with the following properties: 
\begin{itemize} 
\item $s^{(l)}$ solves the Dirac equation at each point in $\hat{\Omega}^{(l)}$.
\item $\chi^{(l)} s^{(l)} = s^{(l)}$ at each point on $\hat{\Sigma}^{(l)}$. 
\item $s^{(l)}$ does not vanish identically. 
\end{itemize} 
Standard unique continuation arguments imply that $\int_U |s^{(l)}|^2 \, d\text{\rm vol}_g > 0$ if $l$ is sufficiently large. By scaling, we may arrange that 
\begin{equation} 
\label{normalization}
\int_U |s^{(l)}|^2 \, d\text{\rm vol}_g = m \, \text{\rm vol}_g(U) 
\end{equation} 
if $l$ is sufficiently large. 

Combining Proposition \ref{L2}, Proposition \ref{sobolev.trace.estimate}, and Proposition \ref{fefferman.phong}, we conclude that 
\begin{align} 
\label{fefferman.phong.estimate}
&\int_{\hat{\Sigma}^{(l)}} \max \{\|d\hat{N}^{(l)}\|_{\text{\rm tr}} - H,0\} \, F^2 \, d\sigma_g \notag \\ 
&\leq o(1) \int_{\hat{\Omega}^{(l)}} |\nabla F|^2 \, d\text{\rm vol}_g + o(1) \int_U F^2 \, d\text{\rm vol}_g 
\end{align}
for every smooth function $F: \hat{\Omega}^{(l)} \to \mathbb{R}$. For each $l$, we apply (\ref{fefferman.phong.estimate}) with $F = (\delta^2+|s^{(l)}|^2)^{\frac{1}{2}}$, and send $\delta \to 0$. This gives 
\begin{align} 
\label{consequence.of.fefferman.phong.estimate}
&\frac{1}{2} \int_{\hat{\Sigma}^{(l)}} \max \{\|d\hat{N}^{(l)}\|_{\text{\rm tr}} - H,0\} \, |s^{(l)}|^2 \, d\sigma_g \notag \\ 
&\leq o(1) \int_{\hat{\Omega}^{(l)}} |\nabla s^{(l)}|^2 \, d\text{\rm vol}_g + o(1) \int_U |s^{(l)}|^2 \, d\text{\rm vol}_g.
\end{align} 
On the other hand, Proposition 2.9 in \cite{Brendle} implies 
\begin{align} 
\label{integral.formula}
&\int_{\hat{\Omega}^{(l)}} |\nabla s^{(l)}|^2 \, d\text{\rm vol}_g + \frac{1}{4} \int_{\hat{\Omega}^{(l)}} R \, |s^{(l)}|^2 \, d\text{\rm vol}_g \notag \\ 
&\leq \frac{1}{2} \int_{\hat{\Sigma}^{(l)}} (\|d\hat{N}^{(l)}\|_{\text{\rm tr}} - H) \, |s^{(l)}|^2 \, d\sigma_g \\ 
&\leq \frac{1}{2} \int_{\hat{\Sigma}^{(l)}} \max \{\|d\hat{N}^{(l)}\|_{\text{\rm tr}} - H,0\} \, |s^{(l)}|^2 \, d\sigma_g. \notag 
\end{align} 
By assumption, the scalar curvature is nonnegative. Combining (\ref{normalization}), (\ref{consequence.of.fefferman.phong.estimate}), and (\ref{integral.formula}), we conclude that 
\[\int_{\hat{\Omega}^{(l)}} |\nabla s^{(l)}|^2 \, d\text{\rm vol}_g \to 0\] 
as $l \to \infty$. Consequently, the sequence $s^{(l)} = (s_1^{(l)},\hdots,s_m^{(l)})$ converges in $C_{\text{\rm loc}}^\infty(\Omega \setminus \partial \Omega,\mathcal{E})$ to an $m$-tuple of parallel spinors $s = (s_1,\hdots,s_m)$ which is defined on the interior of $\Omega$. In particular, we can find a fixed matrix $z \in \text{\rm End}(\mathbb{C}^m)$ such that $\langle s_\alpha,s_\beta \rangle = z_{\alpha\beta}$ at each point in the interior of $\Omega$. Arguing as in Section 4 in \cite{Brendle}, we can show that $z$ is a scalar multiple of the identity matrix. Using (\ref{normalization}), we conclude that $z$ is the identity matrix. 

To summarize, $s = (s_1,\hdots,s_m)$ is a collection of parallel spinors which are defined at each point in the interior of $\Omega$ and are orthonormal at each point in the interior of $\Omega$. This implies that the Riemann curvature tensor of $g$ vanishes identically. 

Since $s$ is parallel, $s$ can be extended continuously to $\Omega$. We next consider the behavior of $s$ along the boundary faces. To that end, we recall the following result from \cite{Brendle}.

\begin{lemma}[cf. \cite{Brendle}, Lemma 3.2]
\label{boundary.faces}
Let $j \in \{0,1,\hdots,q\}$. Then the set 
\[\{u_j=0\} \cap \bigcap_{i \in \{0,1,\hdots,q\} \setminus \{j\}} \{u_i<0\}\] 
is non-empty. Moreover, this set is a dense subset of $\Omega \cap \{u_j=0\}$.
\end{lemma} 

Lemma \ref{boundary.faces} follows from Assumption \ref{no.redundant.inequalities} together with the assumption that $\Omega$ has non-empty interior. We refer to \cite{Brendle} for a detailed proof. \\

Let us fix an arbitrary integer $j \in \{0,1,\hdots,q\}$. Arguing as in \cite{Brendle}, we conclude that 
\[\nu_j \cdot s_\alpha = \sum_{a=1}^n \sum_{\beta=1}^m \langle N_j,E_a \rangle \, \omega_{a\alpha\beta} \, s_\beta\] 
at each point in $\{u_j=0\} \cap \bigcap_{i \in \{0,1,\hdots,q\} \setminus \{j\}} \{u_i<0\}$. In view of Lemma \ref{boundary.faces}, the identity 
\[\nu_j \cdot s_\alpha = \sum_{a=1}^n \sum_{\beta=1}^m \langle N_j,E_a \rangle \, \omega_{a\alpha\beta} \, s_\beta\] 
holds at each point in $\Omega \cap \{u_j=0\}$. Therefore, if $0 \leq j < k \leq q$ and $x$ is a point in $\Omega \cap \{u_j=0\} \cap \{u_k=0\}$, then 
\begin{align*} 
m \, \langle \nu_j,\nu_k \rangle 
&= \sum_{\alpha=1}^m \langle \nu_j,\nu_k \rangle \, \langle s_\alpha,s_\alpha \rangle \\ 
&= \frac{1}{2} \sum_{\alpha=1}^m \langle \nu_j \cdot s_\alpha,\nu_k \cdot s_\alpha \rangle + \frac{1}{2} \sum_{\alpha=1}^m \langle \nu_k \cdot s_\alpha,\nu_j \cdot s_\alpha \rangle \\ 
&= \frac{1}{2} \sum_{a,b=1}^n \sum_{\alpha,\beta,\gamma=1}^m \langle N_j,E_a \rangle \, \langle N_k,E_b \rangle \, \omega_{a\alpha\beta} \, \overline{\omega_{b\alpha\gamma}} \, \langle s_\beta,s_\gamma \rangle \\ 
&+ \frac{1}{2} \sum_{a,b=1}^n \sum_{\alpha,\beta,\gamma=1}^m \langle N_j,E_a \rangle \, \langle N_k,E_b \rangle \, \overline{\omega_{a\alpha\beta}} \, \omega_{b\alpha\gamma} \, \langle s_\gamma,s_\beta \rangle \\ 
&= m \, \langle N_j,N_k \rangle
\end{align*} 
at the point $x$. In the last step, we have used the fact that $\omega_a \omega_b + \omega_b \omega_a = -2\delta_{ab} \, \text{\rm id}$ for all $a,b \in \{1,\hdots,n\}$.

Finally, let us fix an integer $j \in \{0,1,\hdots,q\}$. The arguments in Section 4 in \cite{Brendle} show that the second fundamental form of the hypersurface $\{u_j=0\}$ vanishes at each point in $\{u_j=0\} \cap \bigcap_{i \in \{0,1,\hdots,q\} \setminus \{j\}} \{u_i<0\}$. Using Lemma \ref{boundary.faces}, we conclude that the second fundamental form of the hypersurface $\{u_j=0\}$ vanishes at each point in $\Omega \cap \{u_j=0\}$. This completes the proof of Theorem \ref{main.thm}. 

\appendix 

\section{An estimate for the area of a level set of a convex function}

\begin{proposition} 
\label{area.estimate}
Let $m \geq 2$ be an integer, let $f: \mathbb{R}^m \to \mathbb{R}$ be a convex function, and let $\Sigma := \{f=0\} \cap \{\nabla f \neq 0\}$. Then $\Sigma$ is a smooth hypersurface in $\mathbb{R}^m$. Moreover, $|\Sigma \cap B_r(p)| \leq C(m) \, r^{m-1}$.
\end{proposition}

\textbf{Proof.} 
Let us fix a unit vector $\omega \in S^{m-1}$. Let $\Sigma_\omega := \{f=0\} \cap \{\langle \omega,\nabla f \rangle > 0\}$ and let $V_\omega := \{y \in \mathbb{R}^m: \langle \omega,y-p \rangle = 0\}$ denote the hyperplane orthogonal to $\omega$ passing through $p$. We denote by $\pi_\omega: \Sigma_\omega \to V_\omega$ the orthogonal projection from the hypersurface $\Sigma_\omega$ to the hyperplane $V_\omega$. Since $f$ is a convex function, the map $\pi_\omega$ is injective. At each point on $\Sigma_\omega$, the absolute value of the Jacobian determinant of $\pi_\omega$ is given by $|\det D\pi_\omega| = \frac{\langle \omega,\nabla f \rangle}{|\nabla f|}$. Moreover, the image $\pi_\omega(\Sigma_\omega \cap B_r(p))$ is contained in $V_\omega \cap B_r(p)$. Thus, we conclude that 
\[\int_{\Sigma_\omega \cap B_r(p)} \frac{\langle \omega,\nabla f \rangle}{|\nabla f|} \leq |V_\omega \cap B_r(p)| \leq C(m) \, r^{m-1}.\] 
Finally, we integrate over $\omega \in S^{m-1}$. This gives $|\Sigma \cap B_r(p)| \leq C(m) \, r^{m-1}$. \\

\section{Some elementary inequalities}

\begin{lemma} 
\label{elementary.inequality.1}
The function $t \mapsto \frac{t \cos t}{\sin t}$ is monotone decreasing for $t \in (0,\pi)$.
\end{lemma}

\textbf{Proof.} 
Note that $t \geq \sin t \cos t$ for all $t>0$. This implies 
\[\frac{d}{dt} \Big ( \frac{t \cos t}{\sin t} \Big ) = -\frac{t - \sin t \cos t}{\sin^2 t} \leq 0\] 
for all $t \in (0,\pi)$. This completes the proof of Lemma \ref{elementary.inequality.1}. \\

\begin{proposition} 
\label{elementary.inequality.2}
Let $\alpha \in (0,\frac{\pi}{2})$ and $\beta \in [0,\alpha]$. Then 
\[\frac{\sin(2\beta)}{\sin(2\alpha)} - \frac{\sin(2\theta\beta)}{\sin(2\theta\alpha)} \geq 0\] 
for all $\theta \in (0,1)$.
\end{proposition} 

\textbf{Proof.} 
It suffices to prove the assertion for $\beta \in (0,\alpha]$. Using Lemma \ref{elementary.inequality.1}, we obtain 
\[\frac{2\theta\alpha \cos(2\theta\alpha)}{\sin(2\theta\alpha)} - \frac{2\theta\beta \cos(2\theta\beta)}{\sin(2\theta\beta)} \leq 0\] 
for all $\theta \in (0,1)$. This implies 
\[\frac{d}{d\theta} \Big ( \frac{\sin(2\theta\beta)}{\sin(2\theta\alpha)} \Big ) = -\frac{\sin(2\theta\beta)}{\sin(2\theta\alpha)} \, \Big ( \frac{2\alpha \cos(2\theta\alpha)}{\sin(2\theta\alpha)} - \frac{2\beta \cos(2\theta\beta)}{\sin(2\theta\beta)} \Big ) \geq 0\] 
for all $\theta \in (0,1)$. From this, the assertion follows. \\

\begin{proposition} 
\label{elementary.inequality.3}
Let $\alpha \in (0,\frac{\pi}{2})$ and $\beta \in [0,\alpha]$. Then 
\[\Big | \frac{\sin(2\beta)}{\sin(2\alpha)} - \frac{\sin(2\theta\beta)}{\sin(2\theta\alpha)} \Big | \leq \frac{4(\theta-1)\alpha}{\sin(2\alpha) \sin(2\theta\alpha)}\] 
for all $\theta \in [1,\frac{\pi}{2\alpha})$.
\end{proposition} 

\textbf{Proof.} 
Let $\theta \in [1,\frac{\pi}{2\alpha})$. Then $|\sin(2\alpha) - \sin(2\theta\alpha)| \leq 2(\theta-1)\alpha$ and $|\sin(2\beta) - \sin(2\theta\beta)| \leq 2(\theta-1)\beta$. This implies 
\[|\sin(2\theta\alpha) \sin(2\beta) - \sin(2\alpha) \sin(2\theta\beta)| \leq 2(\theta-1)(\alpha+\beta) \leq 4(\theta-1)\alpha.\] 
From this, the assertion follows.


\begin{thebibliography}{9}
\bibitem{Brendle}
S.~Brendle, \textit{Scalar curvature rigidity of convex polytopes,} Invent. Math. 235, 669--708 (2024)

\bibitem{Fefferman-Phong}
C.~Fefferman and D.~Phong, \textit{Lower bounds for Schr\"odinger equations,} Conference on Partial Differential Equations (Saint Jean de Monts, 1982), Conf. No. 7, pp.~1--7, Soc. Math. France, Paris, 1982

\bibitem{Gromov1}
M.~Gromov, \textit{Dirac and Plateau billiards in domains with corners,} Central European Journal of Mathematics 12, 1109--1156 (2014)

\bibitem{Gromov2}
M.~Gromov, \textit{Four Lectures on Scalar Curvature,} arxiv:1908.10612v6

\bibitem{Gromov3}
M.~Gromov, \textit{Convex Polytopes, dihedral angles, mean curvature, and scalar curvature,} arxiv:2207.13346

\bibitem{Li1}
C.~Li, \textit{A polyhedron comparison theorem for $3$-manifolds with positive scalar curvature,} Invent. Math. 219, 1--37 (2020)

\bibitem{Li2}
C.~Li, \textit{The dihedral rigidity conjecture for $n$-prisms,} J. Diff. Geom. 126, 329--361 (2024)

\bibitem{Wang-Xie-Yu}
J.~Wang, Z.~Xie, and G.~Yu, \textit{On Gromov's dihedral extremality and rigidity conjectures,} arxiv:2112.01510
\end{thebibliography}
\end{document}